\documentclass[epsfig,latexsym,amsfonts,twoside]{article}
\usepackage{amssymb,amsmath,amscd,epsfig}

\def\part#1{\frac{\partial\phantom{#1}}{\partial#1}}
\newtheorem{thm}{Theorem}
\newtheorem{prp}[thm]{Proposition}
\newtheorem{lem}[thm]{Lemma}

\newtheorem{cnj}[thm]{Conjecture}
\newenvironment{prf}{\begin{trivlist}\item[]{\bf Proof} }%
{\hfill $\Box$ \end{trivlist}}
\newenvironment{dfn}{\begin{trivlist}\item[]{\bf Definition}\em }%
{\end{trivlist}}
\newenvironment{rmk}{\begin{trivlist}\item[]{\bf Remark} }%
{\end{trivlist}}
\newenvironment{exm}{\begin{trivlist}\item[]{\bf Example} }%
{\end{trivlist}}

\def\Z{\ifmmode{{\mathbb Z}}\else{${\mathbb Z}$}\fi}
\def\Q{\ifmmode{{\mathbb Q}}\else{${\mathbb Q}$}\fi}
\def\C{\ifmmode{{\mathbb C}}\else{${\mathbb C}$}\fi} 
\def\P{\ifmmode{{\mathbb P}}\else{${\mathbb P}$}\fi} 

\def\H{\ifmmode{{\mathrm H}}\else{${\mathrm H}$}\fi} 

\def\B{\ifmmode{{\cal B}}\else{${\cal B}$}\fi} 
\def\E{\ifmmode{{\cal E}}\else{${\cal E}$}\fi} 
\def\F{\ifmmode{{\cal F}}\else{${\cal F}$}\fi} 
\def\K{\ifmmode{{\cal K}}\else{${\cal K}$}\fi} 
\def\L{\ifmmode{{\cal L}}\else{${\cal L}$}\fi} 
\def\M{\ifmmode{{\cal M}}\else{${\cal M}$}\fi} 
\def\N{\ifmmode{{\cal N}}\else{${\cal N}$}\fi} 
\def\O{\ifmmode{{\cal O}}\else{${\cal O}$}\fi} 
\def\U{\ifmmode{{\cal U}}\else{${\cal U}$}\fi}
\def\X{\ifmmode{{\cal X}}\else{${\cal X}$}\fi} 

\def\Br{\ifmmode{{\mathrm{Br}}}\else{${\mathrm{Br}}$}\fi} 
\def\OG{\ifmmode{\widetilde{\cal M}_4}\else{$\widetilde{\cal M}_4$}\fi} 
\def\D{\ifmmode{{\cal{D}}_{\mathrm{coh}}^b}\else{${{\cal{D}}_{\mathrm{coh}}^b}$}\fi}
\def\Shah{\ifmmode{\amalg\hspace*{-3.5pt}\amalg}\else{$\amalg\hspace*{-3.5pt}\amalg$}\fi}

\begin{document}

\title{Foliations on hypersurfaces in holomorphic symplectic
    manifolds} 
\author{Justin Sawon}
\date{December, 2008}
\maketitle

\begin{abstract}
Let $Y$ be a hypersurface in a $2n$-dimensional holomorphic symplectic
manifold $X$. The restriction $\sigma|_Y$ of the holomorphic
symplectic form induces a rank one foliation on $Y$. We investigate
situations where this foliation has compact leaves; in such cases we
obtain a space of leaves $Y/F$ which has dimension $2n-2$ and admits a
holomorphic symplectic form. 
\end{abstract}

\section{Introduction}

The aim of this article is to introduce a new kind of `holomorphic
symplectic reduction', which produces a holomorphic symplectic
$(2n-2)$-fold from a holomorphic symplectic $2n$-fold. Locally it can
be seen as a holomorphic version of the symplectic reduction coming
from an $S^1$-action, though the group action does not globalize.

A summary of the construction is as follows. We start with a
$2n$-dimensional holomorphic symplectic manifold $X$ with holomorphic
two-form $\sigma$. Let $Y\subset X$ be a hypersurface. The kernel of
$\sigma|_Y$ defines a rank one foliation $F\subset TY$ on the smooth
locus of $Y$, which we call the {\em characteristic foliation\/}. If
this foliation has compact leaves, then we can 
construct a space of leaves $Y/F$ which will be Hausdorff. Note that
$Y/F$ has dimension $2n-2$ and it may be singular. Moreover,
$\sigma|_Y$ will descend to a non-degenerate holomorphic two-form on
the smooth part of $Y/F$.

Although our main interest in this article is hypersurfaces, we will
also consider some examples where $Y\subset X$ is a submanifold of
codimension $k>1$. Note that $\sigma|_Y$ must have rank at least
$2(n-k)$ at each point of $Y$. We say that $Y$ is {\em coisotropic\/}
if the rank of $\sigma|_Y$ is constant and equal to $2(n-k)$ at all
points of $Y$. In this case the kernel of $\sigma|_Y$ defines a
characteristic foliation $F\subset TY$ of rank $k$. If this foliation
has compact leaves, the space of leaves $Y/F$ will be a (possibly
singular) holomorphic symplectic manifold of dimension $2n-2k$.

The fundamental question is thus: can we find hypersurfaces and/or
submanifolds $Y\subset X$ whose characteristic foliations $F$
have compact leaves? Given a foliation on a projective variety, a
criterion for algebraicity (and hence compactness) of the leaves was
given in a recent theorem of Kebekus, Sol{\'a} Conde, and 
Toma~\cite{ksct07}, building on results of Miyaoka~\cite{miyaoka87}
and Bogomolov and McQuillan~\cite{bm01}. Let $Y\subset X$ be a smooth
hypersurface (so that $F$ is regular everywhere) or a coisotropic
submanifold. It follows from the theorem of Kebekus et al.\ that if
$Y$ is covered by curves $C$ such that $F|_C$ is ample then every leaf
is algebraic (and hence compact) and rationally connected. We will
explain how this ensures that the space of leaves $Y/F$ is not only
Hausdorff, but also smooth.

One potential application of these ideas is to constructing new
examples of holomorphic symplectic manifolds. We will describe how
a birational model of O'Grady's example in
ten-dimensions~\cite{ogrady99} arises in this way, starting from a
hypersurface in the Hilbert scheme $\mathrm{Hilb}^6S$ of six points on
a K3 surface $S$. It should be noted that the space of leaves $Y/F$
here is (presumably) quite singular and one still has to apply
O'Grady's desingularization procedure to arrive at a smooth birational
model of his space.

Another potential application is to the classification of holomorphic
symplectic manifolds up to deformation, and we consider two main
results in this direction. The goal is to describe $X$ in terms of the
lower-dimensional space $Y/F$, thereby building a classification
inductively on the dimension. When $X$ is four-dimensional the space
of leaves $Y/F$ is a holomorphic symplectic surface; if $Y/F$ is
smooth, it must be a K3 or abelian surface. For example, let $Y$ be
the inverse image of the diagonal $\Delta$ in the Hilbert scheme of
two points on a K3 surface $S$,
$$\mathrm{Hilb}^2S=\mathrm{Blow}_{\Delta}(S\times S/\Z_2).$$
Then every leaf is a rational curve $\P^1$, and $Y/F=\Delta$ is of
course isomorphic to the K3 surface $S$. Conversely,
Nagai~\cite{nagai03} proved that if a holomorphic symplectic four-fold
contains a divisor $Y$ which is `numerically like' the one above (in
particularly, it can be contracted to a surface $S$) then $S$ is a K3
surface and $X\cong\mathrm{Hilb}^2S$. We consider Nagai's theorem in
the context of foliations, adding some new observations.

The second main classification result involves holomorphic symplectic
manifolds which are Lagrangian fibrations by Jacobians of curves. Let
$X$ be the total space of the relative compactified Jacobian of a
family $\mathcal{C}\rightarrow\P^2$ of genus two curves over
$\P^2$. Markushevich~\cite{markushevich96} proved that if $X$ is a
holomorphic symplectic four-fold then $\mathcal{C}$
is a complete linear system of curves in a K3 surface $S$; in other
words, $X$ is the Beauville-Mukai integrable
system~\cite{beauville99}, and in particular, $X$ is a deformation of
$\mathrm{Hilb}^2S$. Markushevich constructed the K3 surface $S$ as a
branched double cover of the dual plane
$(\P^2)^{\vee}$. Alternatively, we will explain how the K3 surface $S$
can be described as the space of leaves of the characteristic
foliation on a hypersurface $Y\subset X$; namely, we take $Y$ to be
the total space of $\mathcal{C}$, embedded in $X$ by the relative
Abel-Jacobi map. The author has used this approach to generalize
Markushevich's Theorem to higher dimensions~\cite{sawon08ii}: roughly
speaking, if the total space $X$ of the relative compactified Jacobian
of a family $\mathcal{C}\rightarrow\P^n$ of genus $n$ curves over
$\P^n$ is a 
holomorphic symplectic manifold, then $\mathcal{C}$ is a complete
linear system of curves in a K3 surface $S$, i.e., $X$ is the
Beauville-Mukai integrable system and hence deformation equivalent to
$\mathrm{Hilb}^nS$. For the proof one again takes $Y$ to be the total
space of $\mathcal{C}$, embedded by the relative Abel-Jacobi map,
though this time $Y$ will be of higher codimension and thus one must
also prove that $Y$ is coisotropic.

The paper is structured as follows. In Section 2 we outline the main
constructions and prove some general results concerning characteristic
foliations. In Section 3 we describe a number of examples. In Sections
4, 5, and 6 we discuss how Nagai's Theorem, Markushevich's Theorem,
and O'Grady's example in ten dimensions, respectively, can be
understood in terms of foliations on hypersurfaces. Finally, in
Section 7 we make comparisons with the analogous situation in real
symplectic geometry; in particular, we formulate a holomorphic
analogue of the Weinstein Conjecture, asserting the existence of
compact leaves for the characteristic foliation on a hypersurface of
contact type.

While this paper was being written a preprint of Hwang and
Oguiso~\cite{ho07} appeared. They use methods similar to ours to study
the singular fibres of a Lagrangian fibration; namely, if
$\pi:X\rightarrow\P^n$ is a Lagrangian fibration, and
$\Delta\subset\P^n$ is the discriminant locus parametrizing singular
fibres, then Hwang and Oguiso study the characteristic foliation on
the hypersurface $Y:=\pi^{-1}(\Delta)$ and prove that the leaves are
either elliptic curves or chains of rational curves. This in turn
tells us something about the structure of the generic singular fibre
$X_t$, where $t$ is a generic point of $\Delta$. Hwang and
Viehweg~\cite{hv08} also consider characteristic foliations on
hypersurfaces of general type; we mention their results in
Subsection~2.2.

The author would like to thank Jun-Muk Hwang, Yujiro Kawamata, Stefan
Kebekus, Manfred Lehn, and Keiji Oguiso for many helpful discussions
on the material presented here. The author is also grateful for the
hospitality of the Max-Planck-Institut f{\"u}r Mathematik (Bonn), the
Institute for Mathematical Sciences (Chinese University of Hong Kong),
and the Korea Institute for Advanced Studies, where these results were
obtained.

\section{Characteristic foliations and spaces of leaves}

\subsection{Hypersurfaces and coisotropic submanifolds}

\begin{dfn}
A holomorphic symplectic manifold is a compact K{\"a}hler manifold $X$
with a non-degenerate holomorphic symplectic form
$\sigma\in\H^0(X,\Omega^2_X)$. In particular, the dimension
$\mathrm{dim}X=2n$ is even. `Non-degenerate' means that
$\sigma^{\wedge n}$ trivializes the canonical bundle
$K_X=\Omega^{2n}_X$. In addition, if $X$ is simply-connected and
$\H^0(X,\Omega^2_X)$ is one-dimensional, i.e., generated by $\sigma$,
then we say that $X$ is irreducible.
\end{dfn}

\begin{rmk}
For a general compact complex manifold $X$, $\H^0(X,\Omega^2_X)$ is a
quotient of the deRham cohomology group $\H^2(X,\C)$; but if $X$ is
K{\"a}hler then it follows from Hodge theory that the projection
splits and $\H^0(X,\Omega^2_X)$ is a direct summand of
$\H^2(X,\C)$. Therefore the holomorphic symplectic form $\sigma$ is
$d$-closed.
\end{rmk}

Let $Y\subset X$ be a complex hypersurface, possibly with
singularities. The restriction $\sigma|_Y\in\H^0(Y,\Omega^2_Y)$ must be
degenerate; therefore the morphism of sheaves
$TY\rightarrow\Omega^1_Y$ given by contracting a vector field with
$\sigma|_Y$ will not be injective. Define a distribution $F\subset TY$
as the kernel of this morphism. Since $\sigma$ is $d$-closed the
distribution is integrable, i.e., if $v$ and $u$ are local sections of
$F$ then
$$\sigma|_Y([v,u],w)=d\sigma|_Y(v,u,w)+v\sigma|_Y(u,w)-u\sigma|_Y(v,w)=0$$
so $[v,u]$ is also a local section of $F$.

\begin{dfn}
The characteristic foliation on $Y$ is the foliation defined by the
distribution above. By abuse of notation, we use $F$ to denote both
the distribution and the foliation.
\end{dfn}

\begin{rmk}
This notion of integrability of the distribution should not be
confused with (algebraic) integrability of the foliation, meaning that
the leaves are algebraic. We will discuss the latter shortly.
\end{rmk}

Strictly speaking, a foliation in differential geometry should have
constant rank, though algebraic geometers often consider foliations
which are not regular. In other words, the rank of $F$ may jump up on
closed subsets of $Y$; the following lemma shows that this can only
happen on the singular locus of $Y$.

\begin{lem}
Over the smooth locus $Y_{\mathrm{sm}}$ of $Y$ the foliation $F$ is
locally free of rank one. 
\end{lem}

\begin{prf}
Since $\sigma$ is non-degenerate on $X$, $\sigma|_Y$ must be of
constant rank $2n-2$ as a skew-symmetric form over smooth points of
$Y$, and thus $F=\mathrm{ker}\sigma|_Y$ will be of constant rank
one. More explicitly, if $Y$ is given locally by $z_1=0$ in the
neighbourhood of a smooth point, then by the holomorphic Darboux
Theorem we can find local coordinates such that
$$\sigma=dz_1\wedge dz_2+dz_3\wedge dz_4+\ldots+dz_{2n-1}\wedge dz_{2n},$$
$$\sigma|_Y=dz_3\wedge dz_4+\ldots+dz_{2n-1}\wedge dz_{2n},$$
and $F\subset TY$ is generated by $\frac{\partial\phantom{z_2}}{\partial z_2}$.
\end{prf}

\begin{lem}
Over the smooth locus $Y_{\mathrm{sm}}$, $F$ is isomorphic to
$K_Y^{-1}$. In particular, if $Y$ is normal then
$c_1(F)=c_1(TY)=-c_1(K_Y)$.
\end{lem}

\begin{prf}
Consider the exact sequence
$$0\rightarrow F\rightarrow TY\rightarrow TY/F\rightarrow 0.$$
Over the smooth locus $Y_{\mathrm{sm}}$, the holomorphic two-form
$\sigma|_Y$ induces a non-degenerate skew-symmetric form on $TY/F$, as
we have quotiented out the null directions. This implies that
$\bigwedge^{2n-2} TY/F$ is trivial and $F\cong
\bigwedge^{2n-1}TY=K_Y^{-1}$ over $Y_{\mathrm{sm}}$.
\end{prf}

Now suppose that $Y\subset X$ is a (smooth) complex submanifold of
codimension $k$. Once again we get an integrable distribution, i.e., a
foliation, $F\subset TY$ defined as the kernel of
$\sigma|_Y$. However, even though $Y$ is smooth, the rank of
$\sigma|_Y$ could vary. The rank of $\sigma$
restricted to $T_yY$ must be even, at least $2(n-k)$, and at most
$$2\lfloor\frac{\mathrm{dim}Y}{2}\rfloor=2\lfloor\frac{2n-k}{2}\rfloor.$$
When $k=1$ these minimum and maximum values coincide, but this is no
longer the case for $k\geq 2$. The rank is also semi-continuous:,
meaning that
$$\{y\in Y|\mathrm{rank}\sigma|_{T_yY}\leq m\}$$
is a closed subset of $Y$ for each $m$.

\begin{dfn}
We say that $Y$ is a coisotropic submanifold if the rank of
$\sigma|_Y$ takes its minimum value $2(n-k)$ at all points of $Y$,
i.e., if $T_yY\subset T_yX$ is a coisotropic subspace for all $y\in
Y$.
\end{dfn}

If $Y$ is a coisotropic submanifold then the kernel $F\subset TY$ of
$\sigma|_Y$ will be locally free of rank $k$ over $Y$. We once again
have an exact sequence
$$0\rightarrow F\rightarrow TY\rightarrow TY/F\rightarrow 0$$
which implies that $\bigwedge^k F\cong K_Y^{-1}$, i.e., $c_1(F)=-c_1(K_Y)$.

\subsection{Compactness and rationality of the leaves}

Perhaps the earliest criteria for algebraicity and rationality of the
leaves of a holomorphic foliation are due to
Miyaoka~\cite{miyaoka87}. His results were further developed by
Bogomolov and McQuillan~\cite{bm01}. Here we quote results of Kebekus,
Sol{\'a} Conde, and Toma~\cite{ksct07}.

\begin{thm}[\cite{ksct07}]
\label{kebekus}
Let $Y$ be a normal projective variety, $C\subset Y_{\mathrm{sm}}$ a
complete curve contained in the smooth locus of $Y$, and $F\subset TY$
a foliation which is regular over $C$. Assume that $F|_C$ is an ample
vector bundle. Then the leaf through any point of $C$ is algebraic,
and the closure of the leaf through the generic point of $C$ is
rationally connected. If in addition $F$ is regular everywhere then
all leaves are rationally connected.
\end{thm}

This kind of result is proved by applying Mori's bend-and-break
argument to the curve $C$, thereby producing rational curves in the
leaves of $F$.

Now suppose that $Y$ is a hypersurface in a holomorphic symplectic
manifold $X$, as in the previous subsection, and let $F$ be the
characteristic foliation on $Y$. Since the canonical bundle $K_X$ of
$X$ is trivial, we have
$$K_Y=\mathcal{O}(Y)|_Y$$
by adjunction. 
Let $C$ be a curve which
lies in the smooth locus of $Y$; then
$$F|_C\cong K_Y^{-1}|_C=\mathcal{O}(-Y)|_C.$$
So if $\mathcal{O}(Y)$ is ample then $F|_C$ will be anti-ample and the
hypotheses of the theorem can never be satisfied. The following
conjecture was suggested by Jun--Muk Hwang.~\footnote{In a recent
  preprint~\cite{hv08}, Hwang and Viehweg proved that the
  characteristic foliation on a smooth hypersurface of general type
  cannot be such that all of its leaves are algebraic. This implies
  the conjecture since if $\mathcal{O}(Y)$ is ample, then
  $K_Y=\O(Y)|_Y$ is ample and $Y$ is of general type. Their proof
  involves a global {\'e}tale version of Reeb stability (Theorem~3.2
  in~\cite{hv08}); in the context of our argument here, this
  essentially says that there is a generically finite cover
  $\tilde{Y}$ of $Y$ such that when we pull back the foliation to
  $\tilde{Y}$ the genus of the leaves becomes constant.}

\begin{cnj}
Let $Y$ be a smooth hypersurface in an irreducible holomorphic
symplectic manifold $X$ such that $\mathcal{O}(Y)$ is ample. Then the
characteristic foliation on $Y$ cannot be such that all of its leaves
are algebraic.
\end{cnj}

\begin{rmk}
Certainly some of the leaves of the foliation could be algebraic. Let
the algebraic curve $C\subset Y$ be a leaf of the foliation $F$. Then
the tangent bundle of $C$
$$TC\cong F|_C\cong K_Y^{-1}|_C=\mathcal{O}(-Y)|_C$$
is anti-ample and $C$ must have genus $g\geq 2$. Suppose that
every leaf is a (smooth) genus $g$ curve, and indeed that $Y$ is a
fibration by genus $g$ curves over the base $B$, i.e., with no
singular fibres. The homotopy long exact sequence of $Y\rightarrow B$
gives
$$\ldots\rightarrow\pi_2(C)\rightarrow\pi_2(Y)\rightarrow\pi_2(B)\rightarrow\pi_1(C)\rightarrow\pi_1(Y)\rightarrow\pi_1(B)\rightarrow
1.$$
Since $\mathcal{O}(Y)$ is ample and the irreducible holomorphic
symplectic manifold $X$ is simply-connected, $Y$ must have trivial
fundamental group by the Lefschetz Hyperplane Theorem; but then
$\pi_1(C)$ will be the cokernel of the map
$\pi_2(Y)\rightarrow\pi_2(B)$ between abelian groups, forcing
$\pi_1(C)$ itself to be abelian. This gives a contradiction.

The statement of the conjecture is of course stronger than this
because even if all the leaves of $F$ are algebraic, we cannot expect
$Y$ to be fibred by genus $g$ curves: the genus could drop on closed
subsets. The local example
$$Y:=(C\times\C^2)/\Z_2$$
illustrates this, where $C$ is a curve admitting a fixed-point free
involution $\tau$, and $\Z_2$ acts as $\tau$ on $C$ and multiplication
by $-1$ on $\C^2$. Here the central fibre is isomorphic to $C/\tau$
and it has smaller genus.
\end{rmk}

Returning to Theorem~\ref{kebekus}, it is clear that we must
start with a hypersurface $Y$ whose corresponding line bundle
$\mathcal{O}(Y)$ is negative in some sense, in order to apply this
result. Let us assume that this is the case, and that the hypotheses
of the theorem are satisfied. If furthermore $Y$ is smooth then $F$ is
regular everywhere, and every leaf is rationally connected. Since the
leaves are one-dimensional and smooth, they must be rational curves
$\P^1$; the space of leaves will then be particularly well-behaved.

\begin{lem}
Suppose $Y$ is a smooth hypersurface such that all of the leaves of
the characteristic foliation $F$ are rational curves. Then the space
of leaves, denoted $Y/F$, is smooth.
\end{lem}

\begin{prf}
Holmann~\cite{holmann80} proved that if all of the leaves of a
holomorphic foliation on a K{\"a}hler manifold are compact, then the
foliation is stable, meaning that every open neighbourhood of a leaf
$L$ contains a neighbourhood consisting of a union of leaves (known as
a saturated neighbourhood). Stability of the foliation is equivalent
to the space of leaves $Y/F$ being Hausdorff.

Let us describe the local structure of the space of leaves $Y/F$. Let
$L$ be a (compact) leaf of $F$, represented by a point in $Y/F$. Take
a small slice $V$ in $Y$ transverse to the foliation, with $L$
intersecting $V$ at a point $0\in V$. The holonomy map is a group
homomorphism from the fundamental group of $L$ to the group of
automorphisms of $V$ which fix $0$. The holonomy group $H(L)$ of $L$
is the image of the holonomy map. Then $V/H(L)$ is a local model for
the space of leaves $Y/F$ in a neighbourhood of the point representing
$L$ (see Holmann~\cite{holmann78} for details). Globally, the space of
leaves $Y/F$ is constructed by patching together these local models.

If all of the leaves are rational curves, then they are simply
connected and hence the holonomy groups must be trivial. The local
models for $Y/F$ are therefore simply the transverse slices $V$; in
particular, they are smooth.
\end{prf}

\begin{rmk}
Note that $Y$ is a $\P^1$-bundle over the space of leaves $Y/F$. This
is really a holomorphic version of the Reeb Stability Theorem, in the
case when the leaves are smooth and simply connected. By contrast, a
real foliation with leaves isomorphic to $S^1$ need not be an
$S^1$-bundle as some leaves could have multiplicity (e.g. a
Seifert-fibred space with exceptional fibres; sometimes this is
regarded as an $S^1$-bundle over an orbifold). Similarly, a
holomorphic foliation with leaves isomorphic to elliptic curves could
have multiple fibres. In a different direction, if we assume the
leaves are simply connected but not necessarily smooth (i.e., if the
foliation is not regular), then chains of $\P^1$s could occur.
\end{rmk}

\begin{rmk}
The space of leaves $Y/F$ could also be defined in terms of the Chow
variety of subvarieties of $Y$, or in terms of the Hilbert scheme of
subschemes of $Y$; in either case, we take the irreducible component
which contains the generic leaf. When all of the leaves are rational
curves, these constructions of $Y/F$ all yield isomorphic spaces.
\end{rmk}

The situation becomes somewhat more complicated when $Y$ is a
hypersurface with singularities. In this case, the foliation $F$ is
regular along the smooth locus $Y_{\mathrm{sm}}$ of $Y$, but not
regular along the singular locus $Y_{\mathrm{sing}}$ of $Y$. Therefore
it is perhaps better to consider the foliation only over
$Y_{\mathrm{sm}}$; a leaf could then be quasi-projective, i.e., it
could be algebraic but we would need to add a point or points from
$Y_{\mathrm{sing}}$ to compactify it.

Defining the space of leaves is also more delicate. One approach is to
think of the foliation as an equivalence relation on $Y$. When $Y$ is
smooth the equivalence classes are simply the leaves, and this leads
to the local models for $Y/F$ described earlier. When $Y$ is singular
we should take closures of leaves, and moreover, if two closures
intersect we should consider them part of the same equivalence
class. Thus an equivalence class could consist of a chain of
$\P^1$s. Such a construction is similar to Mumford's geometric
invariant theory: if a group $G$ acts on a space $X$, then two points
are equivalent if the closures of their orbits under $G$
intersect. One problem with this approach is that a continuous family
of leaves in $Y_{\mathrm{sm}}$ could all be compactified by the same
point from $Y_{\mathrm{sing}}$; this then produces an equivalence
class of dimension greater than one, which may be undesirable.

We could instead define the space of leaves in terms of the Chow
variety or Hilbert scheme. We cannot guarantee that these different
approaches will yield the same space $Y/F$, but when the generic leaf
is a smooth curve then these various spaces should at least be
birational. Certainly we don't expect $Y/F$ to be smooth in general:
it will usually be singular along points parametrizing leaves which
intersect $Y_{\mathrm{sing}}$.

It is somewhat easier to deal with the case of a higher codimension
submanifold $Y$, provided $Y$ is coisotropic so that the
characteristic foliation $F$ is regular. Firstly, the leaves must be
smooth. If the hypotheses of Theorem~\ref{kebekus} are
satisfied, then all the leaves must be algebraic (hence compact) and
rationally connected. It is well known that smooth rationally
connected varieties are simply connected (see Corollary 4.18 of
Debarre~\cite{debarre01}, for example), and therefore we obtain a
smooth space of leaves $Y/F$ just as in the proof of the lemma above
(i.e., all holonomy groups are trivial so the local models must be
smooth).

We conclude this subsection with the following result about the
structure of the space of leaves $Y/F$.

\begin{lem}
Assume that the space of leaves $Y/F$ is smooth (for example, $Y$ is a
smooth hypersurface or a coisotropic submanifold and the hypotheses of
Theorem~\ref{kebekus} are satisfied). Then $Y/F$ admits a holomorphic
symplectic form.
\end{lem}

\begin{prf}
Intuitively, by taking the space of leaves we are quotienting along
the leaves of $F$, which are the null directions of
$\sigma|_Y$. Therefore $\sigma|_Y$ should descend to a non-degenerate
two-form on $Y/F$. To make this rigourous, we will define holomorphic
symplectic forms on the local models for $Y/F$ and then prove that
they agree when the local models are patched together.

A local model for the space of leaves $Y/F$ in a neighbourhood of the
point representing the leaf $L$ is given by a small slice $V_1$ in
$Y$, with $L$ intersecting $V_1$ transversally at a point $0$. In a
neighbourhood of $0$, $V_1$ will be transverse to the foliation $F$
and therefore the restriction
$\sigma|_{V_1}$ will be a non-degenerate two-form. Suppose $V_2$ is a
second small slice transverse to $L$; the intersections points $L\cap
V_1$ and $L\cap V_2$ may or may not coincide. Shrinking $V_1$ and
$V_2$ if necessary, there is an isomorphism $\phi:V_1\rightarrow V_2$
given by taking $L^{\prime}\cap V_1$ to $L^{\prime}\cap V_2$, where
$L^{\prime}$ is an arbitrary leaf intersecting $V_1$ and $V_2$. This
map is well-defined because the leaves are simply connected, and it
takes the point $0$ in $V_1$ (i.e., $L\cap V_1$) to $0$ in $V_2$
(i.e., $L\cap V_2$). This makes explicit the method of patching
together local models; see Figure~1.

\begin{figure}
\begin{center}
\includegraphics[width=40.0mm]{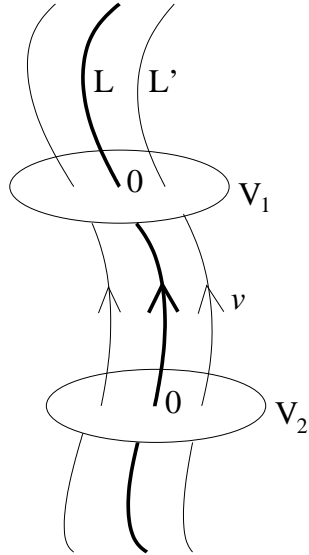}
\end{center}
\caption{Identifying different local models of $Y/F$}
\end{figure}

Now $\phi$ could be regarded as the time $t=1$ map of a flow $\phi_t$
associated to a vector field $v$ along the leaves of the foliation
$F$. Since $d\sigma=0$ and $i(v)(\sigma|_Y)=0$, the Lie derivative
$$\mathcal{L}_v(\sigma|_Y) = v(d\sigma|_Y)-d(i(v)\sigma|_Y)$$
vanishes, and therefore the flow $\phi_t$ will preserve the
holomorphic two-form $\sigma|_Y$. In particular
$$\phi^*(\sigma|_{V_2})=\sigma|_{V_1}$$
which completes the proof.
\end{prf}

\begin{rmk}
In general, the same proof shows that there is always a
holomorphic symplectic form on the smooth part of the space of leaves
$Y/F$. Moreover, if $V/H(L)$ is a local model for $Y/F$ near a leaf
$L$ with non-trivial holonomy group $H(L)$ (in particular, $L$ itself
cannot be simply connected), then $H(L)$ should be a subgroup of the
symplectomorphism group of $V$, i.e., the group of automorphisms
preserving the symplectic structure on $V$. The simplest case would be
when $H(L)$ is a finite group; $V/H(L)$ is then a symplectic
singularity, which may admit a symplectic desingularization. More
generally, foliations can behave very wildly and in principle $H(L)$
could be infinite.
\end{rmk}

\subsection{Divisorial contractions}

Assume now that $X$ is a four-fold, and $Y$ is a hypersurface, not
necessarily smooth or even normal, with $F$ the characteristic
foliation on $Y$. Suppose that the hypotheses of Theorem~\ref{kebekus}
are satisfied, i.e., there is a curve $C$ lying in the smooth locus
$Y_{\mathrm{sm}}$ of $Y$ for which $F|_C$ is ample. Recall that $F$ is
locally free of rank one on $Y_{\mathrm{sm}}$, and
$$F\cong K_Y^{-1}\cong\mathcal{O}(-Y)|_Y$$
over $Y_{\mathrm{sm}}$. Therefore $F|_C$ ample will imply $Y.C<0$, and
the divisor $Y$ on $X$ is not nef.

Let us begin instead with the assumption $q(Y,Y)<0$, where $q$
is the Beauville-Bogomolov quadratic form of $X$. This will imply the
existence of a curve $C$ such that $F|_C$ is ample, as we now
explain. A nef divisor $W$ is the limit of ample divisors, and so it
must satisfy $q(W,W)\geq 0$ since $q(H,H)>0$ for ample $H$. So if $Y$
is not in
$$\{W\mbox{ an effective divisor on }X|q(W,W)\geq 0\}$$
(which contains the nef cone) then $Y$ cannot be nef. It follows that
there must exist a curve $C$ such that $Y.C<0$, or in other words,
such that $F|_C$ is ample. Note that if $Y$ is singular the curve $C$
won't necessarily lie in the smooth locus $Y_{\mathrm{sm}}$, but this
is not important for the following discussion.

We can now run the log minimal model programme (MMP) on
$(X,\epsilon Y)$, where $\epsilon\in\mathbb{Q}_{>0}$. The goal is to reach a
particular birational model for $(X,\epsilon Y)$ after a series of
directed flips. In fact it is known that birational maps between
holomorphic symplectic four-folds factor through Mukai flops, i.e.,
blowing up a $\P^2$ and then blowing it down along a different ruling
of the exceptional divisor (this was proved by Burns, Hu, and
Luo~\cite{bhl03} provided the indeterminacy of the map is normal, but
this hypothesis can be removed because of a later result of Wierzba
and Wi{\'s}niewski~\cite{ww03}). Note that a Mukai flop with respect
to the trivial canonical bundle $K_X$ is a directed flip with respect
to $K_X+\epsilon Y$.

Almost all parts of the log MMP have been proved in dimension
four. The existence of log flips is due to
Shokurov~\cite{shokurov03}. While termination is not yet established
in full generality, there are many partial results: Kawamata, Matsuda,
and Matsuki~\cite{kmm87} proved termination for terminal flips,
Matsuki~\cite{matsuki91} proved termination for terminal flops,
Fujino~\cite{fujino04} proved termination for canonical log flips, and
Alexeev, Hacon, and Kawamata~\cite{ahk07} proved termination for klt
log flips under some additional hypotheses (see also Corti et
al.~\cite{corti07}). In our case $(X,\epsilon Y)$ will be a canonical
pair for sufficiently small $\epsilon$, since $X$ is smooth; we can
therefore apply Fujino's result. The conclusion is that after finitely
many flips we get $(X^{\prime},\epsilon Y^{\prime})$ which, a priori,
satisfies one of the following:
\begin{enumerate}
\item $X^{\prime}$ is a Mori fibre space,
\item $Y^{\prime}$ is nef,
\item $(X^{\prime},\epsilon Y^{\prime})$ admits a divisorial contraction
$$\begin{array}{ccc}
   Y^{\prime} & \subset & X^{\prime} \\
   \downarrow & & \downarrow \\
   S & \subset & \tilde{X}. \\
   \end{array}$$
\end{enumerate}
In the first case $K_{X^{\prime}}+\epsilon Y^{\prime}$ should be
negative on all the fibres, which is impossible since $K_{X^{\prime}}$
is trivial and $Y^{\prime}$ is effective. Next consider case two. By
assumption $q_X(Y,Y)<0$; but then
$$q_{X^{\prime}}(Y^{\prime},Y^{\prime})=q_X(Y,Y)<0$$
since this value is preserved under a Mukai flop. Thus $Y^{\prime}$
cannot be nef either, and case two is impossible. We conclude that the
log MMP must produce case three, a divisorial contraction.

Kaledin~\cite{kaledin06} proved that symplectic resolutions are
semi-small, which implies that $2\mathrm{codim}Y^{\prime}\geq
\mathrm{codim}S$. In other words, $S$ is a surface and the generic
fibre of $Y^{\prime}\rightarrow S$ is one-dimensional; moreover, if
two-dimensional fibres occur they must
be isolated, since $Y^{\prime}$ is irreducible. If all the fibres are
one-dimensional then part (i) of Theorem~1.3 in
Wierzba~\cite{wierzba03} states that the generic fibre must be a tree
of $\P^1$s whose dual graph is a Dynkin diagram; moreover, these
$\P^1$s are (closures of) leaves of the characteristic foliation
$F^{\prime}$ on $Y^{\prime}$. Since $Y^{\prime}$ is irreducible, the
generic fibre must consist of either a single rational curve or a pair
of rational curves joined at a node; we call these type I and type II
fibres respectively. For type II fibres there must be some kind of
monodromy which interchanges the two rational curves, so that
$Y^{\prime}$ is still irreducible; moreover, $Y^{\prime}$ will be
singular along the family of nodes, so if $Y^{\prime}$ is normal this
case cannot occur. Part (ii) of Theorem~1.3 in~\cite{wierzba03} states
that for type I fibres, $S$ is a smooth holomorphic symplectic surface
and $Y^{\prime}\rightarrow S$ is a flat morphism with constant
fibres. Indeed, $S$ is precisely the space of leaves
$Y^{\prime}/F^{\prime}$. A similar statement holds for type II fibres,
though $S$ could now have isolated singularities. In addition, if
there are isolated two-dimensional fibres, the points of $S$ over
which they occur could be singular.

This essentially tells us everything about the characteristic
foliation on the divisor $Y^{\prime}$ in $X^{\prime}$, but we'd really
like to say something about the characteristic foliation on the
original hypersurface $Y$ in $X$. As mentioned
earlier, the birational map $X^{\prime}\dashrightarrow X$ must factor
through a sequence of Mukai flops
$$\begin{array}{ccccccc}
   & \hat{X}_1 & & \ldots & & \hat{X}_m & \\
   & \swarrow\phantom{XX}\searrow & & \swarrow\phantom{XX}\searrow & &
   \swarrow\phantom{XX}\searrow & \\
 X^{\prime}=X_0 & & X_1 & & X_{m-1} & & X_m=X. \\
\end{array}$$
Recall that a Mukai flop blows up a $\P^2$ in $X_k$; the
exceptional divisor in $\hat{X}_{k+1}$ is then isomorphic to 
$\P(\Omega^1_{\P^2})\cong\P(\Omega^1_{(\P^2)^{\vee}})$, and this is
blown down to the dual plane $(\P^2)^{\vee}$ in $X_{k+1}$. Let us
consider the effect of this Mukai flop on the divisor $Y^{\prime}$. If
$Y^{\prime}$ contains (a closure of) a leaf $C_1\cong\P^1$ which
intersects the plane $\P^2$ at a single point, then the Mukai flop
will add a second component $C_2\cong\P^1$ to $C_1$; $C_2$ will lie in
the dual plane $(\P^2)^{\vee}$. A subsequent Mukai flop could later
remove the component $C_2$; or it might remove $C_1$, if it lies in
the plane $\P^2$ which is subsequently flopped.

\begin{exm}
Suppose that $Y^{\prime}$ intersects $\P^2$ in a line, such that each
point on this line lies in precisely one leaf. When we blow up $\P^2$,
we will add a second component to each of these leaves. However, when
we blow down along the other ruling the second components will all end
up intersecting in a single point; this is because a line in $\P^2$
parametrizes a pencil of lines in the dual plane $(\P^2)^{\vee}$, all
meeting at a single point. This is illustrated in Figure~2. Note that
the dual plane $(\P^2)^{\vee}$ is contained in the proper transform of
$Y^{\prime}$. 
\end{exm}

\begin{figure}
\begin{center}
\includegraphics[width=100.0mm]{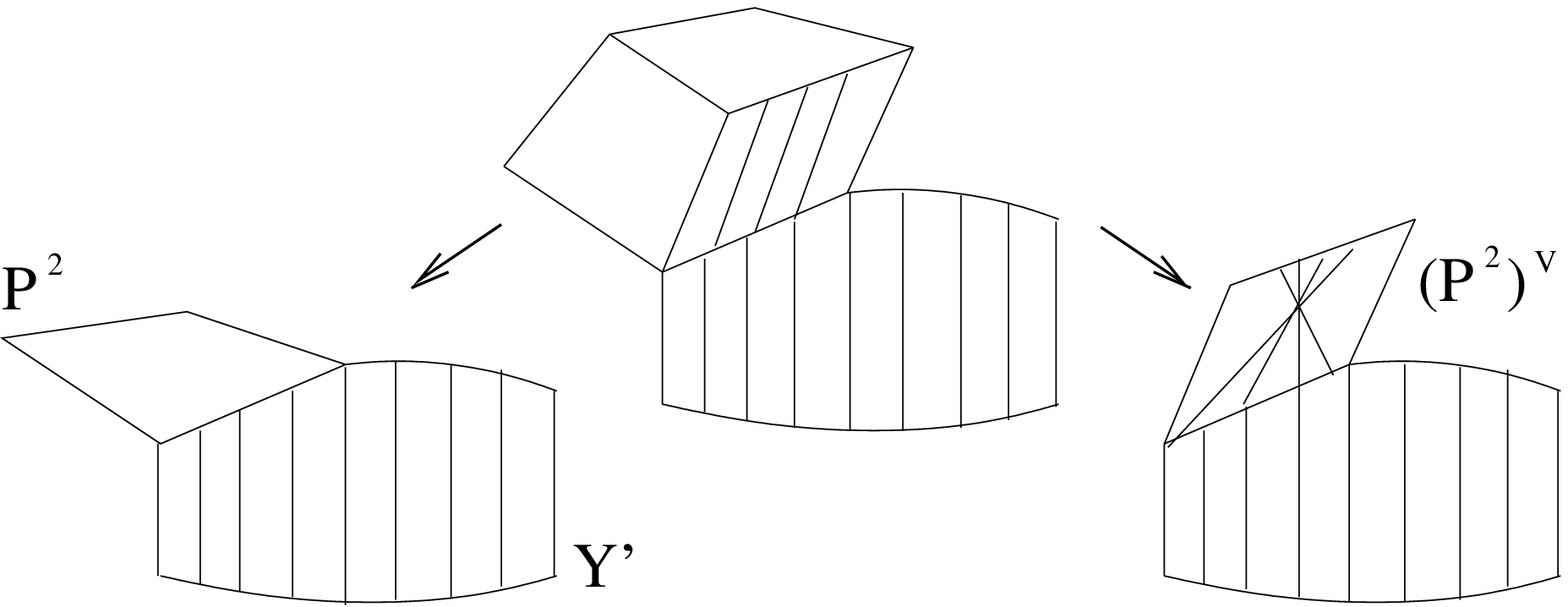}
\end{center}
\caption{Proper transform of $Y^{\prime}$ when $Y^{\prime}\cap\P^2$ is
  a line}
\end{figure}

\begin{exm}
It is also possible for $Y^{\prime}$ to intersect $\P^2$ in a curve
of higher degree, for instance, a conic (an example of this occurs in
the Hilbert scheme of two points on a K3 surface, as described in
Subsection~3.1). Blowing up $\P^2$ again adds
second components to the leaves which intersect $\P^2$ in this
conic. When we blow down along the other ruling the second components
will still intersect each other, but they won't all pass through a
single point. Even more complicated scenarios could arise if the plane
$\P^2$ intersects some leaves in more than one point.
\end{exm}

This gives some idea of how the divisor $Y^{\prime}$ could
change. However, the next lemma shows that the generic leaf is not
disturbed. Denote by $Y_k\subset X_k$ the proper transform of
$Y^{\prime}\subset X^{\prime}$.

\begin{lem}
\label{generic_fibre}
Recall that the generic fibre of $Y^{\prime}\rightarrow S$ is either
of type I (a single rational curves) or of type II (a pair of rational
curves joined at a node), and that these $\P^1$s are (closures of)
leaves of the characteristic foliation on $Y^{\prime}$. The
characteristic foliation on $Y_k$ has generic leaf of the same type;
in particular, this is true for $Y=Y_m$.
It follows that if $Y$ is normal, then $Y_k$ must have
generic leaf of type I for all $k$.
\end{lem}

\begin{prf}
Consider the first Mukai flop $X^{\prime}=X_0\dashrightarrow X_1$. If
$Y^{\prime}$ contained the plane $\P^2$, then we would have a finite
cover $\P^2\rightarrow S$. This contradicts the fact that $S$ is a
holomorphic symplectic surface. Therefore $Y^{\prime}$ and $\P^2$
intersect in a curve, which misses the generic fibre of
$Y^{\prime}\rightarrow S$. Moreover, after the Mukai flop we will
still have a rational map $Y_1\dashrightarrow S$, since the generic
fibre of $Y^{\prime}\rightarrow S$ is not disturbed. This means that
the divisor $Y_1$ cannot contain the plane $\P^2$ which is flopped in
the next birational map $X_1\dashrightarrow X_2$, since otherwise we would
have a rational (multi-valued) map $\P^2\dashrightarrow S$, and a
contradiction as before. Proceeding inductively we see that the
generic fibre of $Y_k\dashrightarrow S$ is never disturbed by the
subsequent Mukai flop $X_k\dashrightarrow X_{k+1}$. This proves the
lemma.
\end{prf}

The next lemma shows that (a closure of) a leaf $C_1$ cannot lie in a
flopped $\P^2$ unless it has previously been replaced by a pair
$C_1+C_2$ of rational curves, as described above.

\begin{lem}
Let $\P^1$ be a generic fibre of $Y_k\dashrightarrow S$ (type I) or
one component of a generic fibre of $Y_k\dashrightarrow S$ (type
II). Then this rational curve $\P^1$ cannot be contained in the plane
$\P^2$ which is subsequently flopped in the birational map
$X_k\dashrightarrow X_{k+1}$.
\end{lem}

\begin{prf}
Since an (analytic) open neighbourhood of this $\P^1$ in $Y_k$ is
isomorphic to an open neighbourhood of the corresponding rational
curve in $Y^{\prime}$, and $Y^{\prime}\rightarrow S$ is a genuine
fibration, we must have
$$N_{\P^1\subset Y_k}\cong\O\oplus\O.$$
Consider the combination of short exact sequences
$$\begin{array}{ccccccccc}
 & & & & 0 & & & & \\
 & & & & \downarrow & & & & \\
 & & & & T\P^1\cong\mathcal{O}(2) & & & & \\
 & & & & \downarrow & & & & \\
 & & & & TX_k|_{\P^1} & & & & \\
 & & & & \downarrow & & & & \\
0 & \rightarrow & N_{\P^1\subset
 Y_k}\cong\mathcal{O}\oplus\mathcal{O} & \rightarrow & N_{\P^1\subset
  X_k} & \rightarrow & N_{Y_k\subset X_k}|_{\P^1} & \rightarrow & 0 \\
 & & & & \downarrow & & & & \\
 & & & & 0 & & & & \\
\end{array}$$
Since $c_1(TX_k)=0$, we must have $N_{Y_k\subset
  X_k}|_{\P^1}\cong\mathcal{O}(-2)$; but then all sequences
  split and
$$TX_k|_{\P^1}\cong\O(-2)\oplus\O\oplus\O\oplus\O(2).$$

On the other hand, if $\P^1\subset\P^2$ then the sequence
$$0\rightarrow T\P^2\rightarrow TX_k|_{\P^2}\rightarrow
N_{\P^2\subset X_k}\cong\Omega^1\P^2\rightarrow 0$$
restricted to $\P^1$ would also split to yield
$$TX_k|_{\P^1}\cong T\P^2|_{\P^1}\oplus\Omega^1\P^2|_{\P^1}$$
which will never look like $\O(-2)\oplus\O\oplus\O\oplus\O(2)$. Note
that we don't need to assume anything about the degree of the rational
curve $\P^1$ in $\P^2$.
\end{prf}

\subsection{Divisors of negative square}

Suppose now that $X$ is a holomorphic symplectic manifold whose
dimension $2n$ could be bigger than four; let $Y\subset X$ be a
hypersurface as before. In higher dimensions some parts of the minimal
model programme remain conjectural at the moment, for example,
termination of flips. However, some of the conclusions of the previous
section are still valid. We will assume that $Y$ has negative square
$q(Y,Y)<0$ with respect to the Beauville-Bogomolov quadratic
form. Note that the quadratic form $q$ is indefinite of signature
$(3,b_2-3)$, where $b_2$ is the second Betti number of $X$. Therefore
whenever $b_2$ is at least four, there will exist classes $\alpha$ in
$\H^2(X,\Z)$ with $q(\alpha,\alpha)<0$; moreover, a standard argument
using the Local Torelli Theorem shows that we can make $\alpha$
algebraic after deforming $X$. A fundamental problem (which we will
not address here) is finding such classes $\alpha$ which can be
represented by effective divisors. The analogous problem for isotropic
classes, i.e., such that $q(\alpha,\alpha)=0$, has been studied, and
effectivity has been established in certain cases (for example,
see~\cite{sawon07}).

By Proposition~5.4 of Huybrechts~\cite{huybrechts03}, if $X$ does not
contain any uniruled divisor then $X$ does not contain any effective
divisors with negative square. Indeed Boucksom later proved that a
prime (i.e., effective and irreducible) divisor with negative square
must be uniruled (Proposition~4.7 of~\cite{boucksom04}). The
uniruledness can be traced back to Huybrechts' arguments, particularly
the fact that the indeterminacy of a birational map $X\dashrightarrow
X^{\prime}$ will be uniruled since the map factors through a blow-up
and blow-down $X\leftarrow Z\rightarrow X^{\prime}$. Here we give an
elementary argument using foliations that yields the same result when
the divisor is smooth.

\begin{prp}
Let $Y$ be a smooth irreducible hypersurface with $q(Y,Y)<0$. Then $Y$
is a $\P^1$-bundle; hence it is uniruled.
\end{prp}

\begin{prf}
Suppose $Y$ were nef. Then $Y$ would be the limit of ample divisors
$H_i$ (at least if we allow rational coefficients). Since
$q(H_i,H_i)>0$ for all $i$, we must have $q(Y,Y)\geq 0$, a
contradiction.

Therefore $Y$ is not nef and there exists a curve $C$, necessarily
contained in $Y$, such that $Y.C<0$. Let $F$ be the characteristic
foliation on $Y$. Then the restriction $F|_C\cong\O(-Y)|_C$ to $C$ is
ample. Since $Y$ is smooth $F$ will be regular everywhere, so all
leaves must be smooth. By Theorem~\ref{kebekus} every leaf of the
foliation is a rational curve $\P^1$, and $Y$ is a $\P^1$-bundle over
the space of leaves $Y/F$.
\end{prf}

If we drop the assumption that $Y$ is smooth then the proof breaks
down, as the curve $C$ may not lie in the smooth locus of $Y$, as
required by Theorem~\ref{kebekus}. Nonetheless, Boucksom's result
implies that $Y$ is still uniruled. The following lemma shows that a
sufficiently 
general rational curve in $Y$ must be (a closure of) a leaf of the
characteristic foliation on $Y$.

\begin{lem}
Let $Y$ be an irreducible hypersurface with $q(Y,Y)<0$; according to
Proposition~4.7 in Boucksom~\cite{boucksom04} $Y$ is uniruled. A
smooth rational curve $\P^1\subset Y$ whose normal bundle
$$N_{\P^1\subset Y}\cong\O(a_1)\oplus\ldots\oplus\O(a_{2n-2})$$
is semi-positive, i.e., $a_i\geq 0$ for all $i$, must be parallel to
the characteristic foliation $F$ on $Y$. Moreover, the normal bundle
must be trivial: $a_i=0$ for all $i$.
\end{lem}

\begin{prf}
Since $TX|_{\P^1}$ admits a non-degenerate symplectic form, it must be
a direct sum of line bundles which are pairwise dual
$$TX|_{\P^1}\cong
\O(-b_1)\oplus\O(b_1)\oplus\ldots\oplus\O(-b_n)\oplus\O(b_n).$$
In particular, $c_1(TX|_{\P^1})=0$ and so the short exact sequence
$$0\rightarrow T\P^1\cong\O(2)\rightarrow TX|_{\P^1}\rightarrow
N_{\P^1\subset X}\rightarrow 0$$
implies that $c_1(N_{\P^1\subset X})=-2$. The exact sequence
$$0\rightarrow N_{\P^1\subset Y}\rightarrow N_{\P^1\subset
  X}\rightarrow N_{Y\subset X}|_{\P^1}\rightarrow 0$$
and the fact that $N_{\P^1\subset Y}$ is a direct sum of line bundles
  $\O(a_i)$ with $a_i\geq 0$ then implies that $c_1(N_{Y\subset
  X}|_{\P^1})\leq -2$. Moreover, this sequence must split so that
$$N_{\P^1\subset X}\cong \O(a_1)\oplus\ldots\oplus\O(a_{2n-2})\oplus
N_{Y\subset X}|_{\P^1}.$$
Substituting this and the decomposition of $TX|_{\P^1}$ into the exact sequence
$$0\rightarrow T\P^1\rightarrow TX|_{\P^1}\rightarrow N_{\P^1\subset
  X}\rightarrow 0,$$
we see that the only possibility is $a_i=0$ for all $i$ and
  $N_{Y\subset X}|_{\P^1}\cong\O(-2)$. Thus $N_{\P^1\subset Y}$ is
  trivial. It follows that the exact sequence
$$0\rightarrow T\P^1\cong\O(2)\rightarrow TY|_{\P^1}\rightarrow
  N_{\P^1\subset Y}\cong\O^{\oplus(2n-2)}\rightarrow 0$$
must also split, and thus
$$F|_{\P^1}\cong K_Y^{-1}|_{\P^1}\cong\bigwedge^{2n-1}TY|_{\P^1}\cong\O(2)$$
and
$$TY/F|_{\P^1}\cong\O^{\oplus(2n-2)}.$$
This means that the diagonal map (the composition of inclusion and
projection) in
$$\begin{array}{ccccccccc}
 & & & & T\P^1 & & & & \\
 & & & & \downarrow & \searrow & & & \\
0 & \rightarrow & F|_{\P^1} & \rightarrow & TY|_{\P^1} & \rightarrow &
 TY/F|_{\P^1} & \rightarrow & 0 \\
\end{array}$$
must be trivial, and hence the inclusion $T\P^1\hookrightarrow
TY|_{\P^1}$ lifts to an isomorphism $T\P^1\cong F|_{\P^1}$. This
completes the proof.
\end{prf}

\section{Examples}

\subsection{The Hilbert square of a K3 surface}

Let $S$ be a K3 surface, and $\Delta\subset S\times S$ be the
diagonal. The symmetric square $\mathrm{Sym}^2S=S\times S/\Z_2$ is
singular along the diagonal, but this can be resolved by
blowing up. The resulting holomorphic symplectic four-fold (discovered
by Fujiki~\cite{fujiki83}) is known as the Hilbert scheme of two
points on $S$, 
$$\mathrm{Hilb}^2S:=\mathrm{Blow}_{\Delta}(S\times S/\Z_2).$$
One obtains the same space by blowing up the diagonal first, and then
quotienting by $\Z_2$.

Let $X$ be this four-fold, and let $Y$ be the inverse image of the
diagonal. Then $Y$ is a $\P^1$-bundle over $\Delta\cong S$. It is easy
to show that these $\P^1$s are precisely the leaves of the
characteristic foliation $F$ on $Y$. Note that in this example $Y$ is
smooth, $F$ is regular, all leaves are $\P^1$s, and the space of
leaves $Y/F\cong S$ is smooth. Note also that as a $\P^1$-bundle over
$S$, $Y$ is the projectivization $\P(TS)$ of the tangent bundle. The
divisor $Y$ is $2-$divisible in the Picard group, i.e., $Y=2\delta$,
but the divisor $\delta$ is not effective; $Y$ is also a divisor of
negative square, $q(Y,Y)=-8$.

\begin{rmk}
If $S$ contains a rational curve, then $\mathrm{Hilb}^2S$ contains
$\mathrm{Sym}^2\P^1\cong\P^2$. The inverse image $Y$ of the diagonal
intersects this in the diagonal of $\mathrm{Sym}^2\P^1$, which is a
conic in the plane $\P^2$, cf.\ the second example in Subsection~2.3.
\end{rmk}

More generally, if $S$ contains a curve $C$ (not necessarily rational) we
could instead choose $Y$ to be the hypersurface 
$$\{\xi\in\mathrm{Hilb}^2S|\mathrm{supp}(\xi)\cap C\not =\emptyset\}.$$
We can divide $Y$ into three strata:
$$Y_1:=\{\xi=\{p,q\}|p\in C, q\in S\backslash C\}\cong C\times
  (S\backslash C),$$
$$Y_2:=\{\xi=\{p,q\}|p,q\in C, p\neq
  q\}\cong\mathrm{Sym}^2C\backslash\Delta,$$
$$Y_3:=\{\xi=\{p,v\}|p\in
  C,v\in\P(T_pS)\cong\P^1\}\cong\P^1\mbox{-bundle over }C.$$
Then $Y_1$ is dense and open, whereas $Y_3$ and the closure of $Y_2$
are divisors in $Y$. Note that $Y$ is not normal: if $\xi=\{p,q\}\in
  Y_2$ then $\xi$ can be deformed in $Y$ by allowing $p$ to vary in
$C$ and $q$ to vary in $S$, or vice versa. This shows that two
branches of $Y$ meet along $Y_2$. 

The leaf through a generic element $\xi=\{p,q\}\in Y_1$ is given by
varying $p$ in $C$ while leaving $q\in S\backslash C$ fixed. The
generic leaf is therefore isomorphic to $C$, and $Y_1$ is a (trivial)
$C$-bundle over $S\backslash C$. To understand the foliation on $Y_2$
and $Y_3$ we really need to take the normalization $\tilde{Y}$ of $Y$,
and pull-back the foliation to $\tilde{Y}$. This mean replacing
$Y_2\cong\mathrm{Sym}^2C\backslash\Delta$ by the double cover $S\times
S\backslash\Delta$; the closure $C\times C$ is then a (trivial)
$C$-bundle over $C$, and the fibres are leaves of the
foliation. By adding $Y_3$ we effectively add a second component,
isomorphic to $\P^1$, to each leaf over $C$. Overall, the space of
leaves on $\tilde{Y}$ is $S$; the generic leaf over a point $q\in
S\backslash C$ is isomorphic to $C$, whereas the `leaf' over $q\in C$ is
isomorphic to the union of $C$ and a rational curve which meets $C$
at the point $q$ (we are regarding the foliation as an equivalence
relation, with each equivalence class being the union of closures of
leaves which intersect, cf. Subsection~2.2). This is illustrated in
Figure~3.

\begin{figure}
\begin{center}
\includegraphics[width=120.0mm]{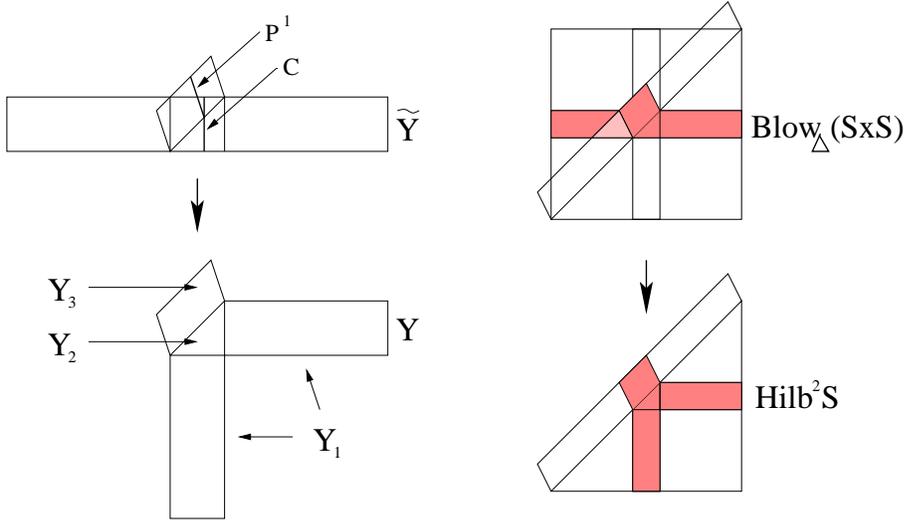}
\end{center}
\caption{Space of leaves on $\tilde{Y}$ is isomorphic to $S$}
\end{figure}

\begin{exm}
When $C$ has genus at least two, this gives an example of a
hypersurface whose characteristic foliation has non-rational
leaves. This does not contradict Conjecture~4 as $Y$ is not smooth
(indeed, not even normal) and moreover $\mathcal{O}(Y)$ is nef but not
ample.
\end{exm}

A higher codimension submanifold of $\mathrm{Hilb}^2S$ is
$$Y:=\{\xi\in\mathrm{Hilb}^2S|p\in\mathrm{supp}(\xi)\}$$
with $p$ some fixed point in $S$. Then $Y$ is isomorphic to $S$
blown up at $p$, and $\sigma|_Y$ is of rank two except on the
exceptional curve where the rank drops to zero ($\sigma|_Y$ vanishes
there). Thus $Y$ is not coisotropic. The only non-trivial leaf is the
exceptional curve and quotienting by the foliation is simply the
blow-down, yielding $S$ as the space of leaves.

\subsection{The Hilbert cube of a K3 surface}

Now suppose that $X$ is the Hilbert scheme $\mathrm{Hilb}^3S$ of three
points on a K3 surface, which is the symplectic desingularization of
$\mathrm{Sym}^3S$ (see Beauville~\cite{beauville83}). Inside
$\mathrm{Sym}^3S$ lies the `big' diagonal $\Delta_{\mathrm{big}}$ where
at least two points coincide and the `small' diagonal
$\Delta_{\mathrm{small}}$ where all three points coincide; thus 
$$\Delta_{\mathrm{big}}:=\{2p+q|p,q\in S\}\cong S\times S$$
and
$$\Delta_{\mathrm{small}}:=\{3p|p\in S\}\cong S.$$
Let $Y$ be the inverse image of $\Delta_{\mathrm{big}}$ in
$\mathrm{Hilb}^3S$. If $p\neq q$, then the fibre of $Y$ above
$2p+q\in\Delta_{\mathrm{big}}$ is precisely $\P(T_pS)\cong\P^1$. These
are the generic leaves of the characteristic foliation on $Y$. On the
other hand, the fibre of $Y$ over
$3p\in\Delta_{\mathrm{small}}\subset\Delta_{\mathrm{big}}$ is the cone
on a rational cubic curve (as described in Subsection~2.5 on punctual
Hilbert schemes in Lehn~\cite{lehn04}). Through a generic point on the
cone the leaf of the foliation is a copy of $\C$, and its closure is a
$\P^1$ obtained by adding the apex of the cone which lies in the
singular locus $Y_{\mathrm{sing}}$ of $Y$. Thus the closures of these
leaves all intersect at the apex of the cone, where the foliation is
not regular. The space of leaves in this example can be defined as 
$$Y/F:=\{\mbox{leaves of }F\}/\sim$$
where two leaves are equivalent if their closures intersect. This
yields $\Delta_{\mathrm{big}}\cong S\times S$ as the space of leaves,
which is a holomorphic symplectic manifold, though not irreducible.

This example can be generalized to higher dimensions. In
$\mathrm{Sym}^nS$ the big diagonal is isomorphic to
$S\times\mathrm{Sym}^{n-2}S$. The characteristic foliation on the
inverse image $Y$ in $\mathrm{Hilb}^nS$ has generic leaf $\P^1$,
though equivalence classes of special leaves become more
complicated. Using the definition above, the space of leaves would be
isomorphic to $S\times\mathrm{Hilb}^{n-2}S$.

\subsection{The generalized Kummer four-fold}

Let $A$ an abelian surface, and $\mathrm{Hilb}^nA$ the Hilbert scheme
of $n$ points on $A$. The composition of the Hilbert-Chow morphism and
addition of points in $A$ gives a map
$$\mathrm{Hilb}^nA\rightarrow\mathrm{Sym}^nA\rightarrow A.$$
The inverse image of $0\in A$ is a holomorphic symplectic manifold
$K_{n-1}$ known as the generalized Kummer
variety~\cite{beauville83}. Once again we take the inverse image of
the big diagonal in $\mathrm{Sym}^nA$, and intersect with $K_{n-1}$ to
get a hypersurface $Y$. As with the earlier examples, the generic leaf
of the characteristic foliation on $Y$ is a $\P^1$.

Consider the case $n=3$, with $Y$ a hypersurface in the four-fold
$K_2$. The point of this example is that even though the generic leaf
is rational, we still have more complicated leaves; in fact there are
$3^4=81$ (equivalence classes of) leaves isomorphic to cones over
cubic rational curves. These occur above the points
$3p\in\Delta_{\mathrm{big}}$ with $p$ a $3$-torsion point in $A$. For
Hilbert schemes of points on K3 surfaces this behaviour only seems to
occur in dimension greater than four.

The generic point of the big diagonal looks like $2p+q$ with $p\neq
q$. If this maps to $0\in A$ under the addition map then $q=-2p$ is
uniquely determined by $p\in A$. The fibre of $Y$ over such a point is
the leaf $\P^1$, and this shows that the space of leaves is the
abelian surface $A$.

\subsection{The Beauville-Mukai system in dimension four}

Let $S$ be a K3 surface containing a smooth genus two curve $C$, and
suppose that $\mathrm{Pic}S$ is generated by $C$. Then $C$ moves in a
two-dimensional linear system $|C|\cong\P^2$. Denote this family of
curves by $\mathcal{C}\rightarrow|C|$, and let 
$$X:=\overline{\mathrm{Pic}}^1(\mathcal{C}/|C|)$$
be the compactified relative Jacobian. Since
$C$ generates $\mathrm{Pic}S$, every curve in the family $\mathcal{C}$
is reduced and irreducible, so the compactified Jacobian is
well-defined (see D'Souza~\cite{dsouza79}). Moreover, $X$ is a smooth
holomorphic symplectic four-fold which is a deformation of
$\mathrm{Hilb}^2S$; it is a Lagrangian fibration over $|C|\cong\P^2$
known as the Beauville-Mukai integrable system~\cite{beauville99}.

We chose the relative Jacobian of degree one so that each curve in
$|C|$ has a canonical embedding into its Jacobian. In other words, the
family of curves $\mathcal{C}$ can be embedded in $X$ by the relative
Abel-Jacobi map; let $Y\subset X$ be the resulting hypersurface.

Write
$$Y\cong\mathcal{C}=\{(C_t,p)|C_t\in|C|,p\in C_t\}.$$
The map $\mathcal{C}\rightarrow|C|$ is given by forgetting $p$, but we
can instead forget $C_t$ which gives a map 
\begin{eqnarray*}
Y & \rightarrow & S \\
(C_t,p) & \mapsto & p. \\
\end{eqnarray*}
The fibre above $p\in S$ is given by the pencil of
curves $C_t$ in $|C|$ which pass through the point $p$; this is a
$\P^1$ for all $p\in S$. Thus $Y$ is a $\P^1$-bundle over $S$. These
$\P^1$s are precisely the leaves of the characteristic foliation $F$
on $Y$. Moreover, $Y$ is again the total space of the projectivization
$\P(TS)$ of the tangent bundle of the K3 surface. To see this note
that $f:S\rightarrow |C|^{\vee}\cong\P^2$ is
a branched double cover of the plane. A direction in the tangent space
$T_pS$ then projects down to a direction in $T_{f(p)}\P^2$, which
defines a line through $f(p)$ in $\P^2$. The inverse image under $f$
of this line is a curve $C_t$ through $p$ in $S$.

Markushevich~\cite{markushevich96} proved the uniqueness of this
example in a certain sense: a holomorphic symplectic four-fold given
by the compactified relative Jacobian of a family of genus two curves
over $\P^2$ (with `mild singularities') must be the Beauville-Mukai
system described above. In Section $5$ we will discuss Markushevich's
Theorem and a generalization due to the author~\cite{sawon08ii}.

Another hypersurface $Y$ in $X\rightarrow |C|$ is given by the inverse
image of a line in $|C|\cong\P^2$. If the line is chosen generically,
$Y$ will be smooth, and hence the characteristic foliation will be
given by $F\cong\mathcal{O}(-Y)|_Y$. Since $\mathcal{O}(Y)$ is nef,
$F$ restricted to a curve $C\subset Y$ can never be ample. Thus
Theorem~\ref{kebekus} will not apply to this example. Moreover, we
claim that the generic leaf is not algebraic; we outline the argument
below.

Let $p\in X$ project down to $t\in\P^2$; in other words, $p$ is in the
fibre $X_t$ of $\pi:X\rightarrow\P^2$. There is an exact sequence
$$0\rightarrow T_pX_t\rightarrow T_pX\rightarrow (N_{X_t\subset
  X})_p\rightarrow 0.$$
The fibre $X_t$ is isomorphic to the Jacobian $\mathrm{Pic}^1(D)$ of
some curve $D$ in the linear system $|C|$, and hence $T_pX_t\cong
\H^0(D,\Omega^1)^{\vee}$. The normal bundle $N_{X_t\subset X}$
corresponds to deformations of $D$ in the K3 surface $S$; since the
normal bundle $N_{D\subset S}\cong\Omega^1_D$, we have the
identification $(N_{X_t\subset X})_p\cong \H^0(D,\Omega^1)$.
The holomorphic symplectic form on $T_pX$ is compatible with the
natural pairing between the dual vector spaces $\H^0(D,\Omega^1)$ and
$\H^0(D,\Omega^1)^{\vee}$.

A line in $|C|\cong\P^2$ corresponds to a pencil of lines through a
fixed point $q^{\prime}$ in the dual plane $|C|^{\vee}$, which in turn 
corresponds to the sublinear system of curves in $|C|$ which all pass
through some fixed point $q$ in the K3 surface $S$. Here $q$ maps to
$q^{\prime}$ under the projection $f:S\rightarrow |C|^{\vee}$. For the
hypersurface $Y$ there is an exact sequence
$$0\rightarrow T_pX_t\rightarrow T_pY\rightarrow (N_{X_t\subset
  Y})_p\rightarrow 0.$$
Here $T_pX_t\cong \H^0(D,\Omega^1)^{\vee}$ as before, but
$(N_{X_t\subset Y})_p\cong \H^0(D,\Omega^1(-q))$ since when we deform
the curve $D$ in $S$ it must still pass through the fixed point $q$.

Restricting the holomorphic symplectic form $\sigma$ to $Y$ gives a
degenerate two-form on $T_pY$. One can show that the kernel lies
entirely in the subspace $T_pX_t\subset T_pY$. Moreover, identifying
$T_pX_t$ with $\H^0(D,\Omega^1)^{\vee}$, we find that
$$\mathrm{ker}\sigma|_Y\cong\mathrm{ker}(\H^0(D,\Omega^1)^{\vee}\twoheadrightarrow\H^0(D,\Omega^1(-q))^{\vee})\cong\H^0(q,\Omega^1|_q)^{\vee}$$
i.e., the foliation $F$ is linear in the abelian surface fibres $X_t$
of $Y$. For a generic curve $D$ through $q$, the vector field in the
direction of $\H^0(q,\Omega^1|_q)^{\vee}$ will be irrational in 
$$X_t\cong\frac{\H^0(D,\Omega^1)^{\vee}}{\H_1(D,\Z)}$$
and hence the leaves of the foliation will be dense in $X_t$, {\em
  not\/} algebraic. Only for some special choices of $D$ will we get
algebraic leaves, which will be isomorphic to elliptic curves.

\subsection{Higher-dimensional Beauville-Mukai systems}

Let $S$ be a K3 surface containing a smooth genus $n$ curve $C$, and
suppose that $\mathrm{Pic}S$ is generated by $C$. Then $C$ moves in an
$n$-dimensional linear system $|C|\cong\P^n$. Denote this family of
curves by $\mathcal{C}\rightarrow|C|$, and let
$$X:=\overline{\mathrm{Pic}}^{n-1}(\mathcal{C}/|C|)$$
be the compactified relative Jacobian. Since $C$ generates
$\mathrm{Pic}S$, every curve in the family $\mathcal{C}$ is reduced
and irreducible, so the compactified Jacobian is well-defined (see
D'Souza~\cite{dsouza79}). Moreover, $X$ is a smooth holomorphic
symplectic manifold of dimension $2n$ which is a deformation of
$\mathrm{Hilb}^nS$; it is a Lagrangian fibration over $|C|\cong\P^n$
known as the Beauville-Mukai system~\cite{beauville99}.

In a smooth fibre $\mathrm{Pic}^{n-1}C$ there is a theta divisor which
can be canonically defined as the image of the map
$\mathrm{Sym}^{n-1}C\rightarrow\mathrm{Pic}^{n-1}C$; it can also be
defined as the locus of degree $n-1$ line bundles $L$ on the curve $C$
which admit at least one non-trivial section, i.e., $h^0(L)>0$. A
generic singular curve in the linear system $|C|$ will have precisely
one node; the compactified Jacobian of such a curve contains a
generalized theta divisor (see Esteves~\cite{esteves97}), which again
can be defined as the locus of rank one torsion-free sheaves which
admit at least one non-trivial section. For a more singular curve
(i.e., worse than simply one node) one could try to define a
generalized theta divisor as above, but it is not clear how
well-behaved it will be. However, these curves with worse
singularities will occur above a codimension two subset
$\Delta_{\mathrm{sing}}\subset |C|$ in the linear system. Therefore a
simpler approach is to take the closure of the relative theta divisor
over $|C|\backslash\Delta_{\mathrm{sing}}$; this produces a
hypersurface $Y\subset X$, and how it looks over
$\Delta_{\mathrm{sing}}$ will usually not be important.

As mentioned above, for a smooth genus $n$ curve the theta divisor in
$\mathrm{Pic}^{n-1}C$ parametrizes line bundles which admit at least
one non-trivial section. A line bundle corresponding to a generic
point in the theta divisor will admit exactly one section up to
scale. Moreover, this section will vanish at $n-1$ distinct points on
$C\subset S$; these points then define a generic element of
$\mathrm{Hilb}^{n-1}S$. These statements also apply to singular curves
in the linear system $|C|$ which have precisely one node. Thus we
obtain a rational map from $Y$ to $\mathrm{Hilb}^{n-1}S$. What is the
generic fibre of this map?

Let $\xi=\{p_1,\ldots,p_{n-1}\}$ be a generic element of
$\mathrm{Hilb}^{n-1}S$ consisting of $n-1$ distinct points. Note that
$S$ is embedded in $|C|^{\vee}\cong\P^n$. The $n-1$ points then
determine a codimension two linear subspace in $|C|^{\vee}$, or
equivalently a pencil of hyperplanes in $|C|^{\vee}$ (the hyperplanes
which contain these $n-1$ points). Each hyperplane intersects $S$ in a
curve which belongs to the linear system $|C|$; the $n-1$ points lie
on this curve and define an effective degree $n-1$ divisor on it. Thus
the pencil of hyperplanes in $|C|^{\vee}$ determines a pencil of
curves in $|C|$, each with a degree $n-1$ line bundle which lies in
the corresponding theta divisor. In other words, the fibre of $Y$
above $\xi\in\mathrm{Hilb}^{n-1}S$ is precisely this pencil, a
rational curve $\P^1$. Finally, one can show that this $\P^1$ is a
leaf of the characteristic foliation $F$ on $Y$.

When $n$ is at least three, the theta divisor in $\mathrm{Pic}^{n-1}C$
will be singular, even for a smooth curve $C$. Therefore we can expect
$Y$ to be singular. Since the characteristic foliation $F$ is not
regular on the singular locus of $Y$, the special leaves will
presumably be fairly complicated.

\subsection{Hurtubise's foliation}

We continue the example from the previous subsection. Instead of a
hypersurface in $X$, we wish to find a coisotropic submanifold of
higher-codimension. Varying the definition slightly, let
$$X:=\overline{\mathrm{Pic}}^1(\mathcal{C}/|C|)$$
be the compactified relative Jacobian of degree one, rather then
degree $n-1$. Every curve, both smooth and singular, can be embedded
in its compactified Jacobian by the Abel-Jacobi map; for the degree
one Jacobian this embedding is canonical, i.e., it does not require
the choice of a basepoint on the curve. Let $Y\subset X$ be the image
of the relative Abel-Jacobi map
$$\begin{array}{ccccc}
 \mathcal{C} & & \hookrightarrow & & X \\
 & \searrow & & \swarrow & \\
 & & |C| & & \\
 \end{array}$$
so that $Y$ is an $(n+1)$-dimensional submanifold of $X$ isomorphic to
the total space of $\mathcal{C}\rightarrow |C|$. Writing
$$Y\cong\mathcal{C}=\{(C_t,p)|C_t\in|C|,p\in C_t\}$$
we see that there is a map given by projection
\begin{eqnarray*}
Y & \rightarrow & S \\
(C_t,p) & \mapsto & p. \\
\end{eqnarray*}
The fibre above $p\in S$ consists of all curves in $|C|$ which pass
through the point $p$. This is a linear codimension one condition on
curves, so the fibre is a hyperplane in $|C|\cong\P^n$, and hence
isomorphic to $\P^{n-1}$ for all $p\in S$. Thus $Y$ is a
$\P^{n-1}$-bundle over $S$; in particular we see that $Y$ is
smooth. Moreover, these $\P^{n-1}$s are precisely the leaves of the
characteristic foliation $F$ on $Y$, and the space of leaves $Y/F$ is
therefore $S$.

Hurtubise considered a local family of Jacobians forming a Lagrangian
fibration in~\cite{hurtubise96}. In other words, he took a family
$\mathcal{C}$ of (smooth) genus $n$ curves over an $n$-dimensional
disc $U\subset\C^n$ such that the total space of the relative Jacobian
$X:=\mathrm{Pic}^d(\mathcal{C}/U)$ admits a holomorphic symplectic
structure, with respect to which the fibres are Lagrangian. We do not
need to compactify since we assume that every curve in the family is
smooth. The degree is also unimportant here, since we can choose a
section of $\mathcal{C}\rightarrow U$; this gives a basepoint in each
curve which allows us to identify Picard groups of different degrees.

Hurtubise then considered the image of the relative Abel-Jacobi map
which embeds the total space of the family $\mathcal{C}\rightarrow U$
as an $(n+1)$-dimensional submanifold in $X$; once again, call the
image $Y\subset X$. Assume that the restriction $\sigma|_Y$ of the
holomorphic symplectic form has constant rank two on $Y$, i.e., assume
that $Y$ is a coisotropic submanifold. The characteristic foliation on
$Y$ then has leaves of dimension $n-1$, and quotienting by the
foliation yields (an open subset, in the analytic topology, of) a
holomorphic symplectic surface $Q$ as the space of leaves
$Y/F$. Hurtubise then shows that $X$ is birational to the Hilbert
scheme $\mathrm{Hilb}^nQ$ of $n$ points on the surface $Q$, with the
rational map identifying the holomorphic symplectic structure on $X$
with the natural holomorphic symplectic structure on
$\mathrm{Hilb}^nQ$ induced from that on $Q$. This result may be viewed
as an example of Sklyanin's separation of variables. Hurtubise and
Markman~\cite{hm98} also proved similar results for fibrations by Prym
varieties.

It is worth pointing out that these are all local results so the issue
of compactness of the leaves of the foliation does not arise. By
studying the compact analogue of this situation, a family of curves
over $\P^n$, and using Theorem~\ref{kebekus} to establish algebraicity
of the leaves of the resulting foliation, the author~\cite{sawon08ii}
was able to generalize Markushevich's Theorem~\cite{markushevich96}
(mentioned in Subsection~3.4 above) to higher dimensions. We will
discuss this result in Section~5.

\subsection{Mukai moduli spaces}

All of the holomorphic symplectic manifolds in the examples above can
be described as Mukai moduli spaces of stable sheaves on K3 or abelian
surfaces~\cite{mukai84}. For example, we can associate to a
zero-dimensional length $n$ subscheme $\xi\in\mathrm{Hilb}^nS$ its
ideal sheaf $\mathcal{I}_{\xi}$; thus $\mathrm{Hilb}^nS$ can be
thought of as the moduli space of rank one stable sheaves with first
Chern class $c_1=0$ and second Chern class $c_2=n$. The hypersurfaces
and coisotropic submanifolds $Y\subset X$ described above can then
often be defined in terms of Brill-Noether loci, e.g. loci of sheaves
with more then the generic number of sections. We first describe the
general setup and then fit some of the above examples into this
framework. We will make a number of assumptions along the way, but the
construction can be modified to accommodate the other cases.

Let $S$ be a K3 surface and let 
$$v=(v_0,v_1,v_2)\in\H^0(S)\oplus\H^2(S)\oplus\H^4(S)$$
be a primitive Mukai vector. Denote by $X:=M^s(v)$ the Mukai moduli
space of stable sheaves $\E$ on $S$ with Mukai vector
$$v(\E):=(r,c_1,r+\frac{c_1^2}{2}-c_2)=v$$
where $r$ is the rank of $\E$, and $c_1$ and $c_2$ are the Chern
classes of $\E$. Since $v$ is primitive, $M^s(v)$ will be a smooth
compact holomorphic symplectic manifold of dimension
$$2+\langle v,v\rangle = 2+v_1^2-2v_0v_2.$$

Assume that the generic sheaf $\E\in M^s(v)$ does not admit any
non-trivial sections. If this were not the case, we could tensor with
an appropriate power of an ample line bundle $L$ on $S$, so that
$\E\otimes L^{-k}$ does not admit any non-trivial sections. Then
instead of the moduli space $M^s(v)$, we could work with the
isomorphic moduli space $M^s(v^{\prime})$ whose elements look like
$\E\otimes L^{-k}$.

Let $Y\subset X=M^s(v)$ be the locus of sheaves $\E$ which admit at
least one non-trivial section, i.e., such that $h^0(\E)>0$. A priori,
$Y$ could be the empty set, but let us assume that this is not the
case (for example, this might be achieved by choosing $k$ above to be
the smallest integer such that $\E\otimes L^{-k}$ does not admit any
non-trivial sections). Moreover, let us assume that the generic
element $\E$ of $Y$ admits exactly one non-trivial section up to
scale, $h^0(\E)=1$ (this seems to be a typical phenomenon), and that
this section gives an injection
$$\O_S\hookrightarrow\E$$
(which is usually the case for generic $\E\in Y$ if the rank $r$ is at
least two).

If all of these assumptions apply, then there is a rational map
$$Y\dashrightarrow M^s(w)$$
given by taking a generic sheaf $\E\in Y$ to the cokernel $\F$ of the
above injection, i.e.,
$$0\rightarrow\O_S\rightarrow\E\rightarrow \F\rightarrow 0.$$
Here the Mukai vector $w$ satisfies
$$w=(w_0,w_1,w_2):=(v_0-1,v_1,v_2-1).$$
The short exact sequence ensures that the cokernel $\F$ has Mukai
vector $w$, though one still has to check that $\F$ will be stable for
generic $\E\in Y$. Note that $M^s(w)$ has dimension
$$2+\langle w,w\rangle = 2+\langle v,v\rangle
-2(1-v_0-v_2)=\mathrm{dim}X-2(1-v_0-v_2).$$
The fibre above a generic element $\F\in M^s(w)$ consists of all $\E$
which can be written as extensions of $\F$ by $\O_S$, up to isomorphism;
but such $\E$ are parametrized by $\P(\mathrm{Ext}^1(\F,\O_S))$, and
thus the fibres are projective spaces. We claim that these projective
spaces are leaves of the characteristic foliation on $Y$.

\begin{prp}
In the above situation, the generic fibre of $Y\dashrightarrow M^s(w)$
is a leaf of the characteristic foliation $F$ on $Y$. Thus the space
of leaves $Y/F$ is birational to $M^s(w)$ and the leaf above a generic
point $\F\in M^s(w)$ is isomorphic to $\P(\mathrm{Ext}^1(\F,\O_S))$.
\end{prp}

\begin{prf}
For $\E\in X=M^s(v)$ the tangent space $T_{\E}X$ can be identified
with $\mathrm{Ext}^1(\E,\E)$. For $\E\in Y$ our goal is to identify the
tangent space $T_{\E}Y\subset T_{\E}X$ to the submanifold
$Y\subset X$ and the distribution $F_{\E}\subset T_{\E}Y$.

We start with the short exact sequence
$$0\rightarrow\O_S\rightarrow\E\rightarrow \F\rightarrow 0$$
and its dual
$$0\rightarrow\F^{\vee}\rightarrow\E^{\vee}\rightarrow\O_S\rightarrow 0$$
(for simplicity, we assume that both $\E$ and $\F$ are locally
free). Tensoring these together and taking the corresponding long
exact sequences gives
$$\begin{array}{ccccccccc}
 & & 0 & & 0 & & 0 & & \\
 & & \downarrow & & \downarrow & & \downarrow & & \\
0 & \rightarrow & {\H}^0(\F^{\vee}) & \rightarrow & \mathrm{Hom}(\F,\E)
& \rightarrow & \mathrm{Hom}(\F,\F) & \rightarrow & \ldots\phantom{a}_{(\clubsuit)} \\
 & & \downarrow & & \downarrow & & \downarrow & & \\
0 & \rightarrow & {\H}^0(\E^{\vee}) & \rightarrow & \mathrm{Hom}(\E,\E)
& \rightarrow & \mathrm{Hom}(\E,\F) & \rightarrow & \ldots\phantom{a}_{(\heartsuit)} \\
 & & \downarrow & & \downarrow & & \downarrow & & \\
0 & \rightarrow & {\H}^0(\O_S) & \rightarrow & {\H}^0(\E)
& \rightarrow & {\H}^0(\F) & \rightarrow & \ldots\phantom{a}_{(\spadesuit)} \\
 & & \downarrow & & \downarrow & & \downarrow & & \\
_{(\clubsuit)}\ldots & \rightarrow & {\H}^1(\F^{\vee}) & \rightarrow & \mathrm{Ext}^1(\F,\E)
& \rightarrow & \mathrm{Ext}^1(\F,\F) & \rightarrow & \ldots\phantom{a}_{(\clubsuit)} \\
 & & \downarrow & & \downarrow & & \downarrow & & \\
_{(\heartsuit)}\ldots & \rightarrow & {\H}^1(\E^{\vee}) & \rightarrow & \mathrm{Ext}^1(\E,\E)
& \rightarrow & \mathrm{Ext}^1(\E,\F) & \rightarrow & \ldots\phantom{a}_{(\heartsuit)} \\
 & & \downarrow & & \downarrow & & \downarrow & & \\
_{(\spadesuit)}\ldots & \rightarrow & {\H}^1(\O_S) & \rightarrow & {\H}^1(\E)
& \rightarrow & {\H}^1(\F) & \rightarrow & \ldots\phantom{a}_{(\spadesuit)} \\
 & & \downarrow & & \downarrow & & \downarrow & & \\
_{(\clubsuit)}\ldots & \rightarrow & {\H}^2(\F^{\vee}) & \rightarrow & \mathrm{Ext}^2(\F,\E)
& \rightarrow & \mathrm{Ext}^2(\F,\F) & \rightarrow & 0 \\
 & & \downarrow & & \downarrow & & \downarrow & & \\
_{(\heartsuit)}\ldots & \rightarrow & {\H}^2(\E^{\vee}) & \rightarrow & \mathrm{Ext}^2(\E,\E)
& \rightarrow & \mathrm{Ext}^2(\E,\F) & \rightarrow & 0 \\
 & & \downarrow & & \downarrow & & \downarrow & & \\
_{(\spadesuit)}\ldots & \rightarrow & {\H}^2(\O_S) & \rightarrow & {\H}^2(\E)
& \rightarrow & {\H}^2(\F) & \rightarrow & 0 \\
 & & \downarrow & & \downarrow & & \downarrow & & \\
 & & 0 & & 0 & & 0 & & \\
\end{array}$$
Note that the first row continues on the fourth row, which continues
on the seventh row; this is indicated by ellipses and $\clubsuit$
(similarly for the other rows).

Since $\E$ and $\F$ are stable, we have
$$\mathrm{Hom}(\F,\F)\cong\C\qquad\mbox{and}\qquad\mathrm{Hom}(\E,\E)\cong\C$$
and then by Serre duality
$$\mathrm{Ext}^2(\F,\F)\cong\C\qquad\mbox{and}\qquad\mathrm{Ext}^2(\E,\E)\cong\C.$$
The slopes of $\E$ and $\F$ must satisfy $0<\mu(\E)<\mu(\F)$ and
therefore
$$\mathrm{Hom}(\F,\E)=0\qquad\mbox{and}\qquad\mathrm{Ext}^2(\E,\F)=0.$$
Since $\mu(\F^{\vee})<\mu(\E^{\vee})<0$ we must have
$${\H}^0(\F^{\vee})=0,\qquad{\H}^0(\E^{\vee})=0,\qquad{\H}^2(\F)=0,\qquad\mbox{and}\qquad{\H}^2(\E)=0.$$
By assumption, a generic $\E\in Y$ admits exactly one non-trivial
section up to scale, and thus
$${\H}^0(\E)\cong\C\qquad\mbox{and}\qquad{\H}^2(\E^{\vee})\cong\C.$$
Since a non-trivial section of $\F$ would lift to a second independent
section of $\E$
$$\begin{array}{ccccccccc}
0 & \rightarrow & {\O}_S & \rightarrow & {\E} & \rightarrow & {\F} &
\rightarrow & 0 \\
& & & & & \nwarrow & \uparrow & & \\
& & & & & & {\O}_S & & \\
\end{array}$$
we must have
$${\H}^0(\F)=0\qquad\mbox{and}\qquad{\H}^2(\F^{\vee})=0.$$
Finally, since $S$ is a K3 surface we have
$${\H}^0(\O_S)\cong\C,\qquad{\H}^1(\O_S)=0,\qquad\mbox{and}\qquad{\H}^2(\O_S)\cong\C.$$

From the third column of the above diagram we can deduce
$$\mathrm{Hom}(\E,\F)\cong\mathrm{Hom}(\F,\F)\cong\C$$
and from the seventh row (or dually) we can deduce
$$\mathrm{Ext}^2(\F,\E)\cong\mathrm{Ext}^2(\F,\F)\cong\C.$$
The diagram now looks like
$$\begin{array}{ccccccccc}
 & & 0 & & 0 & & 0 & & \\
 & & \downarrow & & \downarrow & & \downarrow & & \\
0 & \rightarrow & 0 & \rightarrow & 0 & \rightarrow & {\C} & \rightarrow & \ldots\phantom{a}_{(\clubsuit)} \\
 & & \downarrow & & \downarrow & & \downarrow_{\cong} & & \\
0 & \rightarrow & 0 & \rightarrow & {\C} & \stackrel{\cong}{\rightarrow} & {\C} & \rightarrow & \ldots\phantom{a}_{(\heartsuit)} \\
 & & \downarrow & & \downarrow_{\cong} & & \downarrow & & \\
0 & \rightarrow & {\C} & \stackrel{\cong}{\rightarrow} & {\C}
& \rightarrow & 0 & \rightarrow & \ldots\phantom{a}_{(\spadesuit)} \\
 & & \downarrow & & \downarrow & & \downarrow & & \\
_{(\clubsuit)}\ldots & \rightarrow & {\H}^1(\F^{\vee}) & \rightarrow & \mathrm{Ext}^1(\F,\E)
& \rightarrow & \mathrm{Ext}^1(\F,\F) & \rightarrow & \ldots\phantom{a}_{(\clubsuit)} \\
 & & \downarrow & & \downarrow & & \downarrow & & \\
_{(\heartsuit)}\ldots & \rightarrow & {\H}^1(\E^{\vee}) & \rightarrow & \mathrm{Ext}^1(\E,\E)
& \rightarrow & \mathrm{Ext}^1(\E,\F) & \rightarrow & \ldots\phantom{a}_{(\heartsuit)} \\
 & & \downarrow & & \downarrow & & \downarrow & & \\
_{(\spadesuit)}\ldots & \rightarrow & 0 & \rightarrow & {\H}^1(\E)
& \rightarrow & {\H}^1(\F) & \rightarrow & \ldots\phantom{a}_{(\spadesuit)} \\
 & & \downarrow & & \downarrow & & \downarrow & & \\
_{(\clubsuit)}\ldots & \rightarrow & 0 & \rightarrow & {\C}
& \stackrel{\cong}{\rightarrow} & {\C} & \rightarrow & 0 \\
 & & \downarrow & & \downarrow_{\cong} & & \downarrow & & \\
_{(\heartsuit)}\ldots & \rightarrow & {\C} & \stackrel{\cong}{\rightarrow} & {\C}
& \rightarrow & 0 & \rightarrow & 0 \\
 & & \downarrow_{\cong} & & \downarrow & & \downarrow & & \\
_{(\spadesuit)}\ldots & \rightarrow & {\C} & \rightarrow & 0
& \rightarrow & 0 & \rightarrow & 0 \\
 & & \downarrow & & \downarrow & & \downarrow & & \\
 & & 0 & & 0 & & 0 & & \\
\end{array}$$
from which we extract
$$\begin{array}{ccccccccccc}
 & & & & & & 0 & & & & \\
 & & & & & & \downarrow & & & & \\
0 & \rightarrow & {\C} & \rightarrow & {\H}^1(\F^{\vee}) & \rightarrow
& \mathrm{Ext}^1(\F,\E) & \rightarrow & \mathrm{Ext}^1(\F,\F) & \rightarrow & 0 \\
 & & & & & & \downarrow & & & & \\
 & & & & & & \mathrm{Ext}^1(\E,\E) & & & & \\
 & & & & & & \downarrow & & & & \\
 & & & & & & {\H}^1(\E) & & & & \\
 & & & & & & \downarrow & & & & \\
 & & & & & & 0 & & & & \\
\end{array}$$
Recall that $T_{\E}X\cong\mathrm{Ext}^1(\E,\E)$. An element of the
subspace $\mathrm{Ext}^1(\F,\E)$ will give a deformation of $\E$ that
preserves the short exact sequence description of $\E$, i.e.,
$T_{\E}Y\cong\mathrm{Ext}^1(\F,\E)$. Therefore the vertical sequence
is the normal sequence for $Y\subset X$, and $(N_{Y\subset
  X})_{\E}\cong \H^1(\E)$.

The first column of the larger diagram implies that $\H^1(\E^{\vee})$
is the cokernel of $\C\rightarrow\H^1(\F^{\vee})$. By Serre duality
this is dual to $\H^1(\E)$, and moreover this duality is compatible
with the natural symplectic structure $\sigma$ on
$\mathrm{Ext}^1(\E,\E)$. By exactness an element $u$ of
$$T_{\E}Y\cong\mathrm{Ext}^1(\F,\E)\subset\mathrm{Ext}^1(\E,\E)$$
maps to zero in $\H^1(\E)$, and will therefore pair trivially with any
element $u^{\prime}$ of 
$$\H^1(\E^{\vee})\subset\mathrm{Ext}^1(\F,\E)\cong T_{\E}Y.$$
In other words, $\sigma(u,u^{\prime})=0$ and thus $\H^1(\E^{\vee})$ can be
identified with the null distribution $F_{\E}\subset T_{\E}Y$. It
follows that the diagram above is precisely
$$\begin{array}{ccccccccc}
 & & & & 0 & & & & \\
 & & & & \downarrow & & & & \\
0 & \rightarrow & F_{\E} & \rightarrow
& T_{\E}Y & \rightarrow & (TY/F)_{\E} & \rightarrow & 0 \\
 & & & & \downarrow & & & & \\
 & & & & T_{\E}X & & & & \\
 & & & & \downarrow & & & & \\
 & & & & (N_{Y\subset X})_{\E} & & & & \\
 & & & & \downarrow & & & & \\
 & & & & 0 & & & & \\
\end{array}$$
Finally, observe that
$$F_{\E}\cong\H^1(\E^{\vee})\cong\mathrm{cokernel}(\C\rightarrow\H^1(\F^{\vee}))$$
can be identified with the tangent space to the projective space
$\P(\mathrm{Ext}^1(\F,\O_S))$ at $\E$, and
$(TY/F)_{\E}\cong\mathrm{Ext}^1(\F,\F)$ can be identified with the
tangent space to $M^s(w)$ at $\F$, where of course $\F$ is the image
of $\E$ under the projection $Y\dashrightarrow M^s(w)$. In other
words, $M^s(w)$ is birational to the space of leaves and
$\P(\mathrm{Ext}^1(\F,\O_S))$ is the leaf above a generic point $\F\in
M^s(w)$.
\end{prf}

For our first example, consider the Beauville-Mukai integrable system coming
from a genus $n$ curve $C$ in $S$. A sheaf $L$ on a curve $C$ in the
linear system $|C|$ can be thought of as a torsion sheaf $\iota_*L$ on
the K3 surface $S$ itself, where $\iota:C\hookrightarrow S$ is the
embedding of the curve. This allows us to identify the compactified
relative Jacobian $X=\overline{\mathrm{Pic}}^{n-1}(\mathcal{C}/|C|)$
with the Mukai moduli space $M^s(0,[C],0)$. The generic degree $n-1$
line bundle on a genus $n$ curve will not admit a non-trivial section,
and thus nor will a generic element of $M^s(0,[C],0)$. On the other
hand, the relative theta divisor $Y\subset X$ parametrizes line
bundles $L$ (and more generally, rank one torsion-free sheaves) on
curves $C$ which admit at least one non-trivial section
$$\O_C\stackrel{s}{\longrightarrow}L.$$
This means that the corresponding torsion sheaf $\iota_*L$ on $S$ also
admits a section, given by the composition
$$\O_S\rightarrow\iota_*\O_C\stackrel{\iota_*s}{\longrightarrow} \iota_*L.$$
Thus the relative theta divisor $Y$ is a Brill-Noether locus on
$X=M^s(0,[C],0)$.

In this example the map
$$\O_S\rightarrow\iota_*L$$
is certainly not injective. However, it does fit into an exact sequence
$$0\rightarrow\O_S(-C)\rightarrow\O_S\rightarrow\iota_*L\rightarrow\F\rightarrow
0$$
where $\F$ is defined as the cokernel of the section; but as we vary
the curve in $|C|$, the line bundle $\O_S(-C)$ does not change, so we
can ignore this part of the sequence. Thus we do indeed get a rational
map
$$Y\dashrightarrow M^s(w)$$
which takes $\iota_*L$ to $\F$. Note that $w=(0,0,n-1)$ and
$M^s(w)\cong\mathrm{Hilb}^{n-1}S$, i.e., $\F$ will be the structure
sheaf of a zero-dimensional length $n-1$ subscheme of $S$. Generically
this will just be the set of $n-1$ points on $C$ where the section $s$
of $L$ vanishes.

A similar argument applies when
$X=\overline{\mathrm{Pic}}^1(\mathcal{C}/|C|)$, which is isomorphic to
$M^s(0,[C],-(n-2))$. Once again $Y$, the image of the Abel-Jacobi map
which embeds the family of curves $\mathcal{C}$ in $X$, is a
Brill-Noether locus, parametrizing sheaves $L$ on a curve which admit
non-trivial sections, or equivalently, torsion sheaves $\iota_*L$ on
the K3 surface $S$ which admit non-trivial sections. There is a long
exact sequence
$$0\rightarrow\O_S(-C)\rightarrow\O_S\rightarrow\iota_*L\rightarrow\F\rightarrow
0$$
like before, though now $\F$ has Mukai vector $w=(0,0,1)$, meaning
that it is the structure sheaf of a point. Thus we get a rational map
$$Y\dashrightarrow M^s(w)\cong S$$
as expected.

In these examples, the generic leaves are projective spaces. In
general, Brill-Noether loci in Mukai moduli spaces are stratified and
each open strata is a Grassmannian bundle (see Section~5 of
Markman~\cite{markman01}). The Grassmannian fibres of the largest
strata are the generic leaves of the characteristic foliation on the
Brill-Noether locus; note that Grassmannians are rationally connected,
indicating compatibility with Theorem~\ref{kebekus}. If
$X=\mathrm{Hilb}^nS$ is a Hilbert scheme of points on a K3 surface
$S$, and $Y$ is the inverse image of a (small or big) diagonal in
$\mathrm{Sym}^nS$, then it is less obvious how to describe $Y$ as a
Brill-Noether locus: there does not appear to be a way to {\em
  homologically\/} distinguish a subscheme $\xi\in Y$ from a generic
subscheme $\xi\in X$. However, there may be an isomorphism of
$\mathrm{Hilb}^nS$ with another moduli space on $S$ which takes $Y$ to
a Brill-Noether locus, as in Example~21 in~\cite{markman01}.

\section{Nagai's Theorem}

A possible framework for classifying holomorphic symplectic four-folds
is as follows. Given an irreducible holomorphic symplectic four-fold
$X$, we first find a hypersurface $Y$ in $X$ which should be negative
as a divisor in some sense (certainly $Y$ should not be ample or
nef). In general we will have to deform $X$ in order to find such
$Y$; indeed the generic holomorphic symplectic manifold is
non-projective and contains no hypersurfaces at all. Then we
investigate the characteristic foliation on $Y$ and try to show that
there is a nice space of leaves $S=Y/F$ which is a holomorphic
symplectic surface. Finally, we describe $X$ in terms of the surface
$S$.

Nagai's Theorem~\cite{nagai03} exemplifies this approach. Recall that
inside the Hilbert scheme $X=\mathrm{Hilb}^2S$ of two points on a K3
surface $S$, the inverse image of the diagonal is a $\P^1$-bundle over
the diagonal $\Delta\cong S$. Thus we can recover the K3 surface $S$
as the space of leaves $Y/F$, and clearly $X$ can be described in
terms of $S$. Although he does not explicitly use the language of
foliations, Nagai proved that if a holomorphic symplectic four-fold
$X$ contains a divisor $Y$ satisfying certain hypotheses, then the
space of leaves of the characteristic foliation on $Y$ is a K3 surface
$S$ and $X$ is isomorphic to $\mathrm{Hilb}^2S$. We will describe the
main points of Nagai's proof, adding some new observations. First
recall that by semi-smallness, a divisorial contraction of holomorphic
symplectic four-fold must contract the exceptional divisor to a surface.

\begin{thm}[\cite{nagai03}]
Let $X$ be a projective irreducible holomorphic symplectic
four-fold. Assume there is a birational morphism $f:X\rightarrow Z$
which contracts an irreducible divisor $Y$ to a surface $S\subset Z$
such that
\begin{enumerate}
\item $f|_Y:Y\rightarrow S$ is equidimensional with generic fibre
  $\P^1$,
\item $Y$ is $2$-divisible in $\mathrm{Pic}X$,
\item $Y^4=192.$
\end{enumerate}
Then $S$ is a K3 surface and $X\cong\mathrm{Hilb}^2S$. 
\end{thm}

\begin{rmk}
Nagai notes that if it is known that $X$ is a deformation of the
Hilbert scheme of two points on a K3 surface, then the third condition
may be replaced by $q_X(Y,Y)=-8$, where $q_X$ is the
Beauville-Bogomolov form. We claim that part of the first condition
(that the generic fibre is a rational curve) will also follow from
this assumption.
\end{rmk}

\begin{lem}
Suppose that $X$ is a deformation of the Hilbert scheme of two points
on a K3 surface, and $Y\subset X$ is an irreducible divisor. If
$q_X(Y,Y)=-8$ and $Y$ is $2$-divisible in $\mathrm{Pic}X$, then the
generic leaf of the characteristic foliation on $Y$ is a single
rational curve $\P^1$.
\end{lem} 

\begin{prf}
As described in Subsection~2.3, applying the log minimal model programme to
$(X,\epsilon Y)$ produces, after a sequence of directed flips,
$(X^{\prime},\epsilon Y^{\prime})$ which admits a divisorial
contraction, with the proper transform $Y^{\prime}$ of $Y$ contracted
to a surface. By Lemma~\ref{generic_fibre} the
generic leaves of $Y$ and $Y^{\prime}$ are isomorphic. Now
Wierzba~\cite{wierzba03} proved that the generic fibre of the
contraction of $Y^{\prime}$ to a surface must be a single rational
curve (type I) or a pair of rational curves joined at a node (type
II); suppose it is the latter, and let $C_1$ and $C_2$ be a pair of
rational curves which form a generic fibre. Now a calculation shows
that
$$\O(Y^{\prime})|_{C_1}\cong\O(-1)$$
which contradicts the fact that $Y^{\prime}$ must be $2$-divisible,
like $Y$. Therefore the generic fibre is of type I, completing the
proof.
\end{prf}

\begin{rmk}
Note that we still need to assume that $f|_Y:Y\rightarrow S$ is
equidimensional in order to conclude that $Y$ is a genuine
$\P^1$-bundle over $S$, a fact which is used in several of the
calculations in Nagai's proof. Without this hypothesis, it is possible
that some fibres could have larger dimension, as in the example of the
generalized Kummer four-fold described in Subsection~3.3.
\end{rmk}

Returning to Nagai's Theorem, let us give an outline of his
proof. Since $Y$ is $2$-divisible in $\mathrm{Pic}X$, let $\tilde{X}$
be the double cover of $X$ defined by $\O(\frac{1}{2}Y)$. Let
$\tilde{X}\rightarrow\tilde{Z}\rightarrow Z$ be the Stein
factorization of the composition $\tilde{X}\rightarrow X\rightarrow
Z$. Since $\tilde{X}\rightarrow X$ is ramified over $Y$,
$\tilde{Z}\rightarrow Z$ will be ramified over $f(Y)=S$, around which
$Z$ looks locally like $\C^2\times (A_1\mbox{ surface singularity})$;
it follows that $\tilde{Z}$ will be smooth. Moreover, $\tilde{Z}$ will
have trivial canonical bundle. Some calculations involving the
holomorphic Lefschetz formula of Atiyah and Singer show that
$\tilde{Z}$ is simply connected and
$h^0(\tilde{Z},\Omega^2_{\tilde{Z}})=2$. The Bogomolov Decomposition
Theorem then implies that $\tilde{Z}$ is isomorphic to the product
$S_1\times S_2$ of two K3 surfaces. Finally, Nagai shows that both
$S_1$ and $S_2$ are isomorphic to $S$, which is therefore also a K3
surface, and the covering involution of $\tilde{Z}\rightarrow Z$
simply interchanges the two factors of $\tilde{Z}\cong S\times
S$. This means that we can identify the two diagrams
$$\begin{array}{ccc}
	\tilde{X} & \longrightarrow & \tilde{Z} \\
	\downarrow & & \downarrow \\
	X & \stackrel{f}{\longrightarrow} & Z \\
	\end{array}
	\qquad\qquad\mbox{and}\qquad\qquad 
	\begin{array}{ccc}
	\mathrm{Blow}_{\Delta}(S\times S) & \rightarrow & S\times S \\
	\downarrow & & \downarrow \\
	\mathrm{Hilb}^2S & \rightarrow & \mathrm{Sym}^2S \\
	\end{array}$$
completing the proof.

Note that most of the proof is concerned with reconstructing $X$ from
the space of leaves $S\cong Y/F$; indeed the existence of
the surface $S$ is already part of the hypotheses. In general,
starting with just a hypersurface $Y\subset X$, the difficulty with
constructing a surface $S$ as the image of $Y$ under a divisorial
contraction is that a birational transformation $X\dashrightarrow
X^{\prime}$ might first be required. Of course, one could try to
classify $X^{\prime}$ instead, and note that $X$ is at least
deformation equivalent to $X^{\prime}$, since they are birational (as
proved by Huybrechts~\cite{huybrechts97}). Another approach would be
to show that $X$ does not contain any projective planes $\P^2$; then
no Mukai flop would be possible, so any birational map
$X\dashrightarrow X^{\prime}$ would in fact be an isomorphism. In the
next section we adopt a different approach, though for a different
four-fold $X$ and divisor $Y$: we show that $Y$ is smooth, which will
yield a nice space of leaves $S\cong Y/F$ without the need of any
birational map.

\section{Classification of fibrations by Jacobians}

\begin{dfn}
Let $X$ be an irreducible holomorphic symplectic manifold of dimension
$2n$. We say that $X$ is a Lagrangian fibration if it is fibred over
$\P^n$ with generic fibre an abelian variety of dimension
$n$. Moreover, every fibre should be Lagrangian with respect to the
holomorphic symplectic form.
\end{dfn}

This is essentially the only non-trivial fibration which can exist on
an irreducible holomorphic symplectic manifold, due to results of
Matsushita~\cite{matsushita99, matsushita00} and
Hwang~\cite{hwang07}. Our aim in this section is to consider
Lagrangian fibrations whose fibres are (compactified) Jacobians of
genus $n$ curves. We will recall a classification theorem of
Markushevich~\cite{markushevich96} in dimension four, and describe how
it can be proved using foliations. Then we will describe a
generalization of this classification result to higher dimensions, due
to the author~\cite{sawon08ii}. To begin, we need some control over
the severity of the singular fibres in the fibration.

\begin{dfn}(Markushevich~\cite{markushevich95})
Let $\mathcal{C}\rightarrow\P^2$ be a flat family of genus two
curves. We say that the family has mild degenerations if
\begin{enumerate}
\item the total space $\mathcal{C}$ is smooth,
\item every curve in the family is reduced and irreducible, and the
  singular curves have only nodes or cusps as singularities,
\item if $C_t$ is a singular curve with two singular points $P_1\neq
  P_2$, then the two analytic germs $(\Delta_1,t)$ and $(\Delta_2,t)$
  of the discriminant curves of the unfoldings of the singularities
  $(C_t,P_i)$ meet transversely at $t\in\P^2$.
\end{enumerate}
\end{dfn} 

\begin{rmk}
There is a hypersurface (i.e., curve) $\Delta\subset\P^2$
parametrizing singular curves. Suppose that
$\mathcal{C}\rightarrow\P^2$ has mild degenerations; then $C_t$ will
have a single node if $t$ is a generic point of $\Delta$, it will have
two nodes if $t$ is a node of $\Delta$, and it will have a single cusp
if $t$ is a cusp of $\Delta$. The singularities of $\Delta$ will not
be any worse than this, and nor will the singularities of the curves
$C_t$ be any worse.
\end{rmk}

Since all of the curves in the family are reduced and irreducible
their compactified Jacobians are well-defined.

\begin{thm}(Markushevich~\cite{markushevich95, markushevich96})
Let $\mathcal{C}\rightarrow\P^2$ be a family of genus two curves with
mild degenerations. If $X=\overline{\mathrm{Pic}}^d(\mathcal{C}/\P^2)$
is a holomorphic symplectic four-fold then it is a Beauville-Mukai
integrable system, i.e., the family of curves is a complete linear
system of curves in a K3 surface $S$. In particular, $X$ is a
deformation of the Hilbert scheme $\mathrm{Hilb}^2S$.
\end{thm}

\begin{rmk}
The theorem stated in~\cite{markushevich96} is a strengthening of the
theorem in~\cite{markushevich95}: the degree $d$ is allowed to be
arbitrary instead of just $d=0$, the mild degenerations hypothesis is
relaxed slightly to a necessary and sufficient condition for the
smoothness of the compactified Jacobian
$\overline{\mathrm{Pic}}^d(\mathcal{C}/\P^2)$, and the fact that the
base is isomorphic to $\P^2$ is moved from a hypothesis to a
conclusion. Note that Matsushita~\cite{matsushita99} later proved that
the base of a Lagrangian fibration in four dimensions must always be
isomorphic to $\P^2$. 
\end{rmk}

Let us give a rough outline of Markushevich's proof, which does not
involve foliations. The curves in the family
$\mathcal{C}\rightarrow\P^2$ are hyperelliptic, so the entire family
can be thought of as a branched double cover of a $\P^1$-bundle over
$\P^2$. Moreover, in order for $X$ to be a holomorphic symplectic
four-fold, this $\P^1$-bundle must be the projectivization
$\P(\Omega^1_{\P^2})$ of the cotangent bundle. The branch locus is a
divisor in $\P(\Omega^1_{\P^2})$ which must intersect each $\P^1$
fibre in exactly six points (counted with multiplicity). Again using
the fact that $X$ is a holomorphic symplectic four-fold, one can show
that the branch locus must correspond to a section of
$$\O_{\P(\Omega^1_{\P^2})}(6)\otimes h^*\O(-6)$$
where $h:\P(\Omega^1_{\P^2})\rightarrow\P^2$ is projection. Now
$\P(\Omega^1_{\P^2})$ is the incidence subvariety in
$\P^2\times(\P^2)^{\vee}$ and the line bundle above is actually the
pull-back of $\O_{(\P^2)^{\vee}}(6)$ by the projection to
$(\P^2)^{\vee}$. This means that the branch locus correspond to a
sextic in $(\P^2)^{\vee}$ and one can furthermore show that this
sextic must be smooth. Finally, one obtains the K3 surface $S$ as the
double cover of $(\P^2)^{\vee}$ branched over this sextic;
$\mathcal{C}$ is then a $\P^1$-bundle over $S$ and each curve $C_t$ in
the family projects isomorphically to its image in the K3 surface
$S$.
$$\begin{array}{ccccc}
 & & \mathcal{C} & \rightarrow & S \\
 & \swarrow & \phantom{X}\downarrow{^{_{2:1}}} & & \phantom{X}\downarrow{^{_{2:1}}} \\
{\P}^2 & \leftarrow & {\P}(\Omega^1_{{\P}^2}) & \rightarrow & ({\P}^2)^{\vee} \\
 \end{array}$$
This completes the proof.

Next we outline a different proof of the above theorem due to the
author; the full details may be found in~\cite{sawon08ii}. This proof
uses foliations and works when the degree is one, i.e., when
$X=\overline{\mathrm{Pic}}^1(\mathcal{C}/\P^2)$.

\begin{prf}
The relative Abel-Jacobi map allows us to embed each curve in its
(compactified) Jacobian; when the degree of the Jacobian is one this
embedding is canonical, i.e., does not require the choice of a
basepoint, and therefore the total space of
$\mathcal{C}\rightarrow\P^2$ can be embedded in $X$. Call the
resulting hypersurface $Y\subset X$. Because of the mild degenerations
hypothesis, $Y\cong\mathcal{C}$ is smooth (indeed, this is still the
case even with the slightly weaker conditions on the family of curves
adopted in~\cite{markushevich96}).

Our goal is to show that the leaves of the characteristic foliation
$F$ on $Y$ are smooth rational curves. Then there is a well-defined
space of leaves $S=Y/F$, which we will show is a K3 surface. Moreover,
the curves $C_t$ in the family $\mathcal{C}\cong Y$ will project
isomorphically to their images in $S$, completing the proof. First we
need to find rational curves in $Y$.

Abusing notation, we use $\pi$ to denote the projections of both $Y$
and $X$ to $\P^2$. Let $\ell$ be a generic line in the base $\P^2$,
and let $Z$ be the inverse image $\pi^{-1}(\ell)$ inside $Y$.
$$\begin{array}{ccccc}
Z:=\pi^{-1}(\ell) & \subset & Y\cong\mathcal{C} & \hookrightarrow & X \\
\downarrow & & \downarrow & \swarrow & \\
\ell & \subset & {\P}^2 & & \\
\end{array}$$
Note that $Z$ is a smooth surface which is fibred by genus two curves
over $\ell\cong\P^1$. We claim that $Z$ is {\em not\/} a minimal
surface, i.e., that $Z$ contains $(-1)$-curves. This involves proving
the following statements:
\begin{itemize}
\item $R^1\pi_*\O_X\cong\Omega^1_{\P^2}$ (this was proved by
  Matsushita~\cite{matsushita05}),
\item $R^1\pi_*\O_Y\cong R^1\pi_*\O_X$ (the isomorphism is clear over
  $\P^2\backslash\Delta$, and can be extended over the discriminant
  locus $\Delta$),
\item $R^1\pi_*\O_Z\cong\O_{\ell}(-1)\oplus\O_{\ell}(-2)$ (from
  restricting the previous isomorphisms to $\ell\subset\P^2$),
\item $h^{0,0}(Z)=1$, $h^{0,1}(Z)=0$, and $h^{0,2}(Z)=1$ (from
  substituting the previous isomorphism into the Leray spectral
  sequence for $\pi:Z\rightarrow\ell$).
\end{itemize}
Therefore if $Z$ were a minimal surface, it's Kodaira dimension would
be least zero. However, we can also prove the following:
\begin{itemize}
\item if $X$ is an irreducible holomorphic symplectic four-fold then
  the characteristic number $\sqrt{\hat{A}}[X]$ is at least $25/32$
  (this follows from Guan's bounds~\cite{guan01} on the characteristic
  numbers of holomorphic symplectic four-folds),
\item $\mathrm{deg}\Delta=24(2\sqrt{\hat{A}}[X])^{1/2}$ (as proved by
  the author in~\cite{sawon08i}),
\item $\mathrm{deg}\Delta\geq 30$ and hence $Z\rightarrow\ell$ has at
  least $30$ singular fibres (this follows from the previous two
  points), 
\item $c_2(Z)=\mathrm{deg}\Delta-4\geq 26$ and $K_Z^2\leq -2$ (follows
  from the previous inequality and Noether's formula).
\end{itemize}
Thus if $Z$ were a minimal surface, it's Kodaira dimension would be
$-\infty$. We conclude that $Z$ is not minimal, establishing the
claim.

We have shown that for each $\ell\subset\P^2$ the surface
$Z=\pi^{-1}(\ell)$ contains at least one $(-1)$-curve. Next we show
that these $(-1)$-curves are leaves of the foliation on $Y$. Denoting
one of these curves by $C\cong\P^1$, consider the following
combination of short exact sequences:
$$\begin{array}{ccccccccc}
 & & & & 0 & & & & \\
 & & & & \downarrow & & & & \\
 & & & & TC & & & & \\
 & & & & \downarrow & & & & \\
0 & \rightarrow & F|_C & \rightarrow & TY|_C & \rightarrow & TY/F|_C &
\rightarrow & 0 \\
 & & & & \downarrow & & & & \\
0 & \rightarrow & N_{C\subset Z} & \rightarrow & N_{C\subset Y} &
\rightarrow & N_{Z\subset Y}|_C & \rightarrow & 0 \\
 & & & & \downarrow & & & & \\
 & & & & 0 & & & & \\
 \end{array}$$
We can identify $TC\cong\O(2)$, $N_{C\subset Z}\cong\O(-1)$ since $C$
is a $(-1)$-curve in $Z$, $N_{Z\subset Y}\cong\O(k)$ where $k$ is the
degree of the projection map $\pi|_C:C\rightarrow\ell$, and
$TY/F|_C\cong\O(a_1)\oplus\O(a_2)$ where $a_1+a_2=0$. Substituting
these into the exact sequence, one deduces that we must have
$F|_C\cong\O(2)$ and the map $TC\rightarrow TY|_C$ must lift to an
isomorphism $TC\cong F|_C$. In other words, $C$ is a leaf of the
characteristic foliation on $Y$. It also follows that $k=1$, so that
$C$ is a section of $Z\rightarrow\ell$.

Thus we obtain a space of leaves $S=Y/F$ which is a smooth surface
admitting a holomorphic symplectic structure, i.e., a K3 or abelian
surface. A leaf of the foliation will intersect each curve $C_t$ in
$\mathcal{C}\cong Y$ at most once since the leaf is a section of
$Z\rightarrow\ell$; moreover, one can show that an intersection must
be transverse. Therefore each curve $C_t$ maps isomorphically to its
image in $S$ under the projection $Y\rightarrow S$. Since $S$ contains
a $2$-dimensional linear system of genus two curves, it is a K3
surface, and this completes the proof.
\end{prf}

The advantage of this proof is that it can be generalized to higher
dimensions. First we need an appropriate generalization of a family of
curves with mild degenerations.

\begin{dfn}
Let $\mathcal{C}\rightarrow\P^n$ be a flat family of genus $n$
curves. If $C_t$ is a singular curve of this family with isolated
singular points $P_1,\ldots ,P_k$ then each singular point has a
versal deformation space $\mathcal{X}(P_i)$ and there is an induced
map
$$\phi_t:T_t\P^n\rightarrow T_0\mathcal{X}(P_1)\times\ldots\times
T_0\mathcal{X}(P_k).$$
We say that the family has mild singularities if 
\begin{enumerate}
\item every curve in the family is reduced and irreducible,
\item $\phi_t$ is surjective for every singular curve $C_t$ in the family.
\end{enumerate}
\end{dfn} 

\begin{rmk}
It follows from these conditions that the total space of
$\mathcal{C}\rightarrow\P^n$ is smooth (see~\cite{sawon08ii}). In
fact, smoothness of $\mathcal{C}$ is a sufficient hypothesis for the
following theorem.
\end{rmk}

\begin{thm}(Sawon~\cite{sawon08ii})
\label{main}
Let $\mathcal{C}\rightarrow\P^n$ be a family of genus $n$ curves with
mild singularities. If $X=\overline{\mathrm{Pic}}^1(\mathcal{C}/\P^n)$
is a holomorphic symplectic manifold of dimension $2n$ and the degree
of the discriminant locus $\Delta\subset\P^n$ is greater than $4n+20$
then $X$ is a Beauville-Mukai integrable system, i.e., the family of
curves is a complete linear system of curves in a K3 surface $S$. In
particular, $X$ is a deformation of the Hilbert scheme
$\mathrm{Hilb}^nS$.
\end{thm}

\begin{prf}
The main ideas of the proof are the same as the one above, though some
steps are a little more complicated. We outline them here; full
details may be found in~\cite{sawon08ii}.

As before, the relative Abel-Jacobi embedding
$\mathcal{C}\hookrightarrow\overline{\mathrm{Pic}}^1(\mathcal{C}/\P^n)$
gives a submanifold $Y\subset X$ of dimension $n+1$. Although $Y$ is
smooth, we don't know a priori whether it is coisotropic. However, we
will show later that the rank of the restriction $\sigma|_Y$ is
constant, so that the characteristic foliation $F$ is regular and $Y$
is coisotropic.

Again we look for rational curves in $Y$. Let $\ell$ be a generic line
in the base $\P^n$ and let $Z:=\pi^{-1}(\ell)$. As before, we can show
that $h^{0,0}(Z)=1$, $h^{0,1}(Z)=0$, and $h^{0,2}(Z)=1$. We also have
$c_2(Z)=\mathrm{deg}\Delta-4n+4$ and
$K_Z^2=-\mathrm{deg}\Delta+4n+20< 0$; this is where we use the
bound on the degree of the discriminant locus. As before, we conclude
that $Z$ is {\em not\/} minimal, i.e., it contains at least one
$(-1)$-curve.

Again we analyze the combination of short exact sequences:
$$\begin{array}{ccccccccc}
 & & & & 0 & & & & \\
 & & & & \downarrow & & & & \\
 & & & & TC & & & & \\
 & & & & \downarrow & & & & \\
0 & \rightarrow & F|_C & \rightarrow & TY|_C & \rightarrow & TY/F|_C &
\rightarrow & 0 \\
 & & & & \downarrow & & & & \\
0 & \rightarrow & N_{C\subset Z} & \rightarrow & N_{C\subset Y} &
\rightarrow & N_{Z\subset Y}|_C & \rightarrow & 0 \\
 & & & & \downarrow & & & & \\
 & & & & 0 & & & & \\
 \end{array}$$
We can show that for a generic
$(-1)$-curve $C$, the restriction $F|_C$ will be locally
free. Then analyzing the above sequences as before reveals that
$$F|_C\cong\O(1)^{\oplus(n-2)}\oplus \O(2).$$
In particular $F$ has rank $n-1$ at a generic
point of $Y$. Equivalently, $\sigma|_Y$ has rank two at a generic
point; but then it must have rank two everywhere since the rank can
only drop on closed subsets and two is the minimal value. Therefore
$F$ is locally free, and $Y$ is coisotropic.

To complete the proof, we show that each $(-1)$-curve $C$ is a section
of $Z\rightarrow\ell$. Of course $C$ lies in a leaf of the
foliation $F$. Moreover, by considering the normal bundle of $C$
inside the leaf, we can show that each leaf maps isomorphically to a
hyperplane in $\P^n$ under the projection
$Y\cong\mathcal{C}\rightarrow\P^n$. So the leaves are $\P^{n-1}$s, and
$Y$ is a $\P^{n-1}$-bundle over the space of leaves $S=Y/F$. As
before, $S$ is a smooth surface admitting a holomorphic symplectic
form. The curves $C_t$ in the family $\mathcal{C}\cong Y$ map
isomorphically to their images under the projection $Y\rightarrow
S$. Since $S$ contains an $n$-dimensional linear system of genus $n$
curves, it is a K3 surface, and this completes the proof.
\end{prf}

\section{O'Grady's example in dimension ten}

In~\cite{ogrady99} O'Grady constructed a new ten-dimensional
holomorphic symplectic manifold by desingularizing a certain moduli
space of semi-stable sheaves on a K3 surface $S$. In this section we
show that for certain choices of the K3 surface $S$, O'Grady's example
is birational to the space of leaves of the characteristic foliation
on a hypersurface in the twelve-dimensional holomorphic symplectic
manifold $\mathrm{Hilb}^6S$. The author learnt of this construction
from Nigel Hitchin; apparently it is due to Dominic Joyce.

\begin{thm}[O'Grady~\cite{ogrady99}]
Let $S$ be a K3 surface and let $M^s(2,0,-2)$ be the non-compact Mukai
moduli space of stable rank two sheaves on $S$ with Chern classes
$c_1=0$ and $c_2=4$. Let $M^{ss}(2,0,-2)$ be the compactification
given by adding ($S$-equivalence classes of) semi-stable sheaves; note
that $M^{ss}(2,0,-2)$ is singular along the locus of strictly
semi-stable sheaves. Then $M^{ss}(2,0,-2)$ can be desingularized by a
sequence of blow-ups and blow-downs, such that the resulting space
$\tilde{M}^{ss}(2,0,-2)$ is a smooth holomorphic symplectic manifold
of dimension ten.
\end{thm}

We start by describing how O'Grady's example for a special choice of
K3 surface is birational to a Lagrangian fibration. 

\begin{prp}[O'Grady~\cite{ogrady99}, Rapagnetta~\cite{rapagnetta04}]
Let $S$ be a K3 surface which is a double cover of the plane, branched
over a generic sextic. Then O'Grady's space $\tilde{M}^{ss}(2,0,-2)$
is birational to a Lagrangian fibration.
\end{prp}

\begin{prf}
The birational map to a Lagrangian fibration was first discovered by
O'Grady himself; Rapagnetta later showed that we can assume the
Lagrangian fibration is smooth. Let $H$ be the polarization of $S$
given by the pull-back of a line in the plane; since $S$ is a generic
branched cover of the plane, we can assume that $\mathrm{Pic}S$ is
generated by $H$. Consider the linear system $|2H|\cong\P^5$ on $S$,
whose generic element is a smooth genus five curve. Denote this family
of curves by $\mathcal{C}\rightarrow |2H|$. Denote by
$\overline{\mathrm{Pic}}^6(\mathcal{C}/|2H|)$ the degree six
compactified relative Jacobian of $\mathrm{C}/|2H|$. Since $|2H|$
contains both reducible curves (the pull-back of a pair of lines in
the plane) and non-reduced curves (the pull-back of a double line), we
{\em define\/} $\overline{\mathrm{Pic}}^6(\mathcal{C}/|2H|)$ to be the
irreducible component of the Mukai moduli space of semi-stable sheaves
on $S$ which contains $\iota_*L$, where $\iota:C\hookrightarrow S$ is
the inclusion of a generic curve of the linear system into $S$ and $L$
is a degree six line bundle on $C$. In other words
$$\overline{\mathrm{Pic}}^6(\mathcal{C}/|2H|):=M^{ss}(0,2H,2).$$
Note that this space is singular as the Mukai vector $(0,2H,2)$ is
{\em not\/} primitive.

It can be shown that for generic $C$ and $L$, $h^0(C,L)=2$ and $L$ is
globally generated. Then we can define a sheaf $\E$ as the kernel of
the evaluation map
$$\mathrm{H}^0(C,L)\otimes\O_S\rightarrow \iota_*L$$
and it turns out that $\F:=\E\otimes\O(H)$ is a stable sheaf with
Mukai vector $(2,0,-2)$. This induces a birational map
$$\overline{\mathrm{Pic}}^6(\mathcal{C}/|2H|)\dashrightarrow\tilde{M}^{ss}(2,0,-2)$$
where the space on the left is a Lagrangian fibration over
$|2H|\cong\P^5$. Moreover, Rapagnetta~\cite{rapagnetta04} proved that
$\overline{\mathrm{Pic}}^6(\mathcal{C}/|2H|)$ can be desingularized by
using the same sequence of blow-ups and blow-downs as O'Grady used to
desingularize $M^{ss}(2,0,-2)$. This produces a smooth Lagrangian
fibration which is birational to O'Grady's space.
\end{prf}

We complete this section by using the Lagrangian fibration to show the
following result.

\begin{prp}
Let $S$ be a K3 surface which is a double cover of the plane, branched
over a generic sextic. Then O'Grady's space $\tilde{M}^{ss}(2,0,-2)$
is birational to the space of leaves of the characteristic foliation
on a hypersurface in $\mathrm{Hilb}^6S$.
\end{prp}

\begin{prf}
It suffices to work with the singular space
$\overline{\mathrm{Pic}}^6(\mathcal{C}/|2H|)$ instead of O'Grady's
space $\tilde{M}^{ss}(2,0,-2)$, since they are birational. As
mentioned earlier, if $\iota_*L$ is a generic element of
$\overline{\mathrm{Pic}}^6(\mathcal{C}/|2H|)$, with $L$ a degree six
line bundle on a curve $C\in |2H|$, then $h^0(C,L)=2$. A non-trivial
section of $L$ vanishes at six points of $C\subset S$, counted with
multiplicity, and therefore defines an element of
$\mathrm{Hilb}^6S$. Note that this works even if there is some
repetition among the six points: for example, if $p_1=p_2$ and the
other points are distinct then the length six zero-dimensional subscheme
is given by the four points $p_3,\ldots,p_6\in C\subset S$ plus
$\{p_1,v\}$, where $v\in T_{p_1}S$ is any non-zero vector in
$T_{p_1}C\subset T_{p_1}S$. Moreover, for generic $L$ we won't
encounter vanishing to multiplicity greater than two. Of course, two
sections which agree up to scale will result in the same element of
$\mathrm{Hilb}^6S$. Thus we get a family of elements of
$\mathrm{Hilb}^6S$ parametrized by $\P(\H^0(C,L))\cong\P^1$.

We have shown that an open subset of
$\overline{\mathrm{Pic}}^6(\mathrm{C}/|2H|)$ parametrizes a family of
rational curves in $\mathrm{Hilb}^6S$. The closure of the locus swept
out by this family of curves is a hypersurface $Y$ in
$\mathrm{Hilb}^6S$. This hypersurface can also be defined by
$$Y:=\{\xi\in\mathrm{Hilb}^6S|\exists C\in|2H|\mbox{ such that }\xi\in C\}$$
i.e., $Y$ consists of those length six zero-dimensional subschemes of
$S$ which lie on a curve in the linear system $|2H|$. To see that this
defines a hypersurface, note that if we project six points in $S$ down
to the plane then the condition that they lie on a conic, and hence
that the original six points lie on a curve in $|2H|$, is codimension
one. A generic element $\xi\in Y$ will lie on precisely one curve $C$
in the linear system $|2H|$, and the rational map
$$Y\dashrightarrow\overline{\mathrm{Pic}}^6(\mathcal{C}/|2H|)$$
is given by mapping $\xi$ to the degree six divisor on $C$ which it
defines. The fibres of this map are the rational curves in
$\mathrm{Hilb}^6S$ described earlier, and these are the generic leaves
of the characteristic foliation $F$ on $Y$. Thus we have shown that
$\overline{\mathrm{Pic}}^6(\mathcal{C}/|2H|)$ is birational to the
space of leaves $Y/F$, and this completes the proof.
\end{prf}

\begin{rmk}
It is not clear whether the space of leaves $Y/F$ is smooth, though
this would seem unlikely. Nevertheless, it is possible that it admits
a simpler symplectic desingularization than $M^{ss}(2,0,-2)$; this
would give a more direct route to constructing O'Grady's space.
\end{rmk}

\section{Holomorphic Weinstein conjecture}

In this final section we explain how the existence or non-existence of
compact leaves of the characteristic foliation $F$ on the hypersurface
$Y$ can be regarded as a holomorphic analogue of the Weinstein
Conjecture from (real) symplectic geometry.

First we recall the Weinstein Conjecture~\cite{weinstein79}. Let
$(M,\omega)$ be a (real) symplectic manifold of dimension $2n$, let
$H$ be a smooth function on $M$, and let $X_H$ be the corresponding
Hamiltonian vector field, which is defined by $\omega(X_H,v)=dH(v)$
for $v\in TM$. Since
$$X_H(H)=dH(X_H)=\omega(X_H,X_H)=0$$
the vector field $X_H$ preserves the level sets of $H$. Assume that
$\lambda\in\mathbb{R}$ is chosen so that the level set
$H^{-1}(\lambda)$ is compact. A fundamental question is: will the
vector field $X_H$ have a periodic orbit on $H^{-1}(\lambda)$?

\begin{exm}
Let $M$ be $\mathbb{R}^{2n}$ with the standard symplectic form, and let
$$H=(\mu_1x_1^2+\mu_2x_2^2+\ldots+\mu_{2n}x_{2n}^2)/2$$
with $\mu_1,\ldots,\mu_{2n}$ positive real numbers. Then
$H^{-1}(\lambda)$ is an ellipsoid,
$$dH=\mu_1x_1dx_1+\mu_2x_2dx_2+\ldots+\mu_{2n}x_{2n}dx_{2n}$$
and
$$X_H=\mu_1x_1\frac{\partial\phantom{x_2}}{\partial
  x_2}-\mu_2x_2\frac{\partial\phantom{x_1}}{\partial
  x_1}+\ldots+\mu_{2n-1}x_{2n-1}\frac{\partial\phantom{x_{2n}}}{\partial
  x_{2n}}-\mu_{2n}x_{2n}\frac{\partial\phantom{x_{2n-1}}}{\partial
  x_{2n-1}}.$$
For generic $\mu_1,\ldots,\mu_{2n}$ there will be precisely $n$
  periodic orbits, the intersections of $H^{-1}(\lambda)$ with the
  $\{x_1,x_2\}$-plane, the $\{x_3,x_4\}$-plane, etc.
\end{exm}

The Weinstein Conjecture~\cite{weinstein79} states that there is
always a periodic orbit if $H^{-1}(\lambda)$ is of {\em contact type\/}. The
latter means that $\omega$ restricted to $H^{-1}(\lambda)$ is exact
and equal to $d\alpha$, where $\alpha$ is a contact form on
$H^{-1}(\lambda)$, i.e., $\alpha\wedge (d\alpha)^{n-1}\neq 0$. An
equivalent definition is as follows: $H^{-1}(\lambda)$ is of contact
type if there exists a Liouville vector field $v:U\rightarrow TM|_U$
on an open neighbourhood $U$ of $H^{-1}(\lambda)$ in $M$, i.e., a
vector field $v$ whose Lie derivative preserves the symplectic form,
$\mathcal{L}_v\omega=\omega$. The corresponding flow $\psi_t$ will
rescale the symplectic form by an exponential,
$\psi_t^*(\omega)=e^t\omega$. The Weinstein Conjecture was proved for
$M=\mathbb{R}^{2n}$ by Viterbo~\cite{viterbo87}.

Next we explain how this relates to characteristic foliations on
hypersurfaces. Since $\omega(X_H,-)=dH$, when restricted to the level
set $H^{-1}(\lambda)$ we find that $\omega|_{H^{-1}(\lambda)}(X_H,-)$
is identically zero. In other words, $X_H$ generates the
characteristic foliation on $H^{-1}(\lambda)$; a periodic orbit of
$X_H$ is precisely a compact leaf of the foliation.

In the holomorphic case there are no non-constant functions on a
compact complex manifold $X$. Instead the hypersurface $Y$ can be
regarded as $s^{-1}(0)$ where $s$ is a section of the corresponding
line bundle $\mathcal{O}(Y)$. Locally we can take $z_1=s$ to be the
first coordinate in a system of Darboux coordinates, in which case the
characteristic foliation $F$ is generated locally by
$\frac{\partial\phantom{z_2}}{\partial z_2}$, which is also the
holomorphic analogue of the Hamiltonian vector field corresponding to
$s$. Globally everything needs to be twisted by a power of the line
bundle $\mathcal{O}(Y)$.

The analogue of a periodic orbit isomorphic to $S^1$ is a compact leaf
isomorphic to $\P^1$. As we have seen, the Reeb Stability Theorem
guarantees that nearby leaves will also be rational curves, because
$\P^1$ is simply-connected. This is in contrast to the real case,
where a periodic orbit is not simply-connected and so isolated
periodic orbits can exist. An example is the generic ellipsoid
described above, which has finitely many periodic orbits.

In the holomorphic case, the contact form $\theta$ on the hypersurface
$Y$ also needs to be twisted by a line bundle $L$; thus
$$\theta\in\H^0(Y,\Omega^1_Y\otimes L)$$
and since
$$\theta\wedge(d\theta )^{n-1}\in\H^0(Y,\Omega^{2n-1}_Y\otimes L^n)$$
should be non-vanishing, we must have $K_Y\cong L^{-n}$. More details
on holomorphic contact geometry can be found in~\cite{beauville98},
though note that Beauville considers the Fano case, i.e., when $L$ is
ample. In our case, $L^n\cong K_Y^{-1}\cong\mathcal{O}(-Y)$ should be
semi-positive in some sense, but not necessarily ample.

The following could be called the Holomorphic Weinstein Conjecture.

\begin{cnj}
Let $X$ be a compact holomorphic symplectic manifold and let $Y\subset
X$ be a hypersurface of holomorphic contact type. Then
the generic leaf of the characteristic foliation on $Y$ is a rational
curve. In particular, if $Y$ is smooth then it is a $\P^1$-bundle over
the space of leaves $Y/F$.
\end{cnj}

\begin{rmk}
In our examples, the hypersurface is often observed to be of contact
type. For instance, our hypersurfaces in $\mathrm{Hilb}^2S$ turned out
to be isomorphic to the projectivization $\P(TS)$ of the tangent
bundle of $S$; but since $S$ is a K3 surface,
$TS\cong\Omega^1S$, and it is well known that the canonical
one-form on $\Omega^1S$ induces a contact form on $\P(\Omega^1S)$.
\end{rmk}

\begin{flushleft}
Department of Mathematics\hfill sawon@math.colostate.edu\\
Colorado State University\hfill www.math.colostate.edu/$\sim$sawon\\
Fort Collins CO 80523-1874\\
USA\\
\end{flushleft}

\end{document}